\documentclass[10pt, reqno]{amsart}
\usepackage{graphicx} 
\graphicspath{ {./images/} }
\usepackage[utf8]{inputenc} 
\usepackage{kotex}
\usepackage{dsfont}

\usepackage{amssymb}
\usepackage{amsmath}
\usepackage{mathtools}
\usepackage{tabularx}
\usepackage{enumitem}
\usepackage{setspace}
\usepackage[table,xcdraw]{xcolor}
\usepackage{geometry}
\usepackage{tikz}
\usepackage[colorlinks,linkcolor=purple]{hyperref}
\usepackage{indentfirst}
\usepackage{refcheck}

\usepackage{amsthm}

\newtheorem*{claim}{Claim}

\title[Singular C-S model on general digraphs]{Global existence and emergent dynamics of singular Cucker-Smale model on general digraphs}

\author[Jung]{Jinwook Jung}
\address[Jinwook Jung]{\newline Department of Mathematics and  Research Institute for Natural Sciences \newline
Hanyang University, 222 Wangsimni-ro, Seongdong-gu, Seoul 04763, Republic of Korea}
\email{jinwookjung@hanyang.ac.kr}

\author[Kang]{Hyungsuk Kang}
\address[Hyungsuk Kang]{\newline Department of Mathematics \newline
Hanyang University, 222 Wangsimni-ro, Seongdong-gu, Seoul 04763, Republic of Korea}
\email{sukkun6620@hanyang.ac.kr}

\author[Lim]{Wonjin Lim}
\address[Wonjin Lim]{\newline Department of Mathematics \newline
Hanyang University, 222 Wangsimni-ro, Seongdong-gu, Seoul 04763, Republic of Korea}
\email{robot8305@hanyang.ac.kr}

\author[Shin]{Eojin Shin}
\address[Eojin Shin]{\newline Department of Mathematics \newline
Hanyang University, 222 Wangsimni-ro, Seongdong-gu, Seoul 04763, Republic of Korea}
\email{tlsdjwls2357@hanyang.ac.kr}

\numberwithin{equation}{section}

\newtheorem{theorem}{Theorem}[section]
\newtheorem{lemma}{Lemma}[section]
\newtheorem{corollary}{Corollary}[section]
\newtheorem{proposition}{Proposition}[section]
\newtheorem{remark}{Remark}[section]
\newtheorem{definition}{Definition}[section]

\newcommand{\R}{\mathbb R}

\newcommand{\N}{\mathbb N}

\newcommand{\bq}{\begin{equation}}
\newcommand{\eq}{\end{equation}}

\newcommand{\lt}{\left}
\newcommand{\rt}{\right}

\makeatletter
\def\moverlay{\mathpalette\mov@rlay}
\def\mov@rlay#1#2{\leavevmode\vtop{%
   \baselineskip\z@skip \lineskiplimit-\maxdimen
   \ialign{\hfil$\m@th#1##$\hfil\cr#2\crcr}}}
\newcommand{\charfusion}[3][\mathord]{
    #1{\ifx#1\mathop\vphantom{#2}\fi
        \mathpalette\mov@rlay{#2\cr#3}
      }
    \ifx#1\mathop\expandafter\displaylimits\fi}
\makeatother

\begin{document}

\allowdisplaybreaks

\date{\today}


\begin{abstract}
We investigate the global-in-time dynamics of the Cucker--Smale model under general network topologies with a singular communication weight function $\psi(r) = r^{-\alpha}$, $\alpha>0$. For $\alpha \ge 1$, we establish the global-in-time existence of smooth solutions when network topologies belong to one of several classes of directed graphs. In particular, in one dimension, we present a necessary condition on the network topology for global existence by constructing finite-time blow-up solutions from non-collisional initial configurations. Moreover, we provide a sufficient framework for the emergence of mono-cluster flocking. For certain classes of network topologies, we further show that mono-cluster flocking occurs unconditionally for non-collisional initial configuration when $\alpha =1$. 
\end{abstract}

\maketitle
\tableofcontents

\centerline{\date}

%
%
%
%

\section{Introduction}
\setcounter{equation}{0}

Collective motions in a group of self-propelled agents are ubiquitous in our nature and society, such as aggregation of bacteria, flocking of birds, herding of sheep, etc \cite{Ba08, CHZ17, MT14, To98}. Among them, we are interested in the phenomena called {\it flocking}, where self-driven agents adjust their motion to reach a velocity alignment based on inter-agent communications and limited environmental information. 
The first mathematical modelling of this behavior began in the seminal works by Reynolds and Vicsek et al \cite{R87, V95}. Later, another prototype model for flocking has been introduced by Cucker and Smale \cite{CS07}. Specifically, let $(x_i, v_i) \in \R^d\times\R^d$ be the position and velocity of the $i$-th agent governed by the Cucker-Smale model (in short, CS model), respectively. Then the time evolution of $(x_i, v_i)$ is governed by the following system of ODEs: 

\bq\label{main_sys}
\begin{aligned}
&\dot{x}_i= v_i, \quad i=1, \dots, N, \quad t>0,\\
& \dot{v}_i= \frac{1}{N}\sum_{j\neq i} a_{ij} \psi(|x_j-x_i|)(v_j-v_i)\\
&(x_i(0), v_i(0)) = (x_i^0, v_i^0) \in \R^d \times \R^d,
\end{aligned}
\eq
where $\psi$ is the communication weight function that measures the strength of the interaction between agents based on the inter-agent distances. Here, $A= (a_{ij})$ denotes the 0-1 adjacency matrix of the underlying network topology:

\[
a_{ij} = \left\{\begin{array}{ll} 1 & \mbox{if agent $j$ influences agent $i$,}\\
0 &\mbox{otherwise}.\end{array} \right.
\]
In recent decades, one finds numerous literature on the CS model \eqref{main_sys} in a various point of views, e.g., emergence of flocking with regular communication weight \cite{CS07, HL09, HT08}, collision avoidance \cite{ACHL12, CCMP17, CD10}, control problem \cite{CKPP19, Pe09, So23}, mean-field limits \cite{HKZ18}, general network topologies \cite{CD16, DM11, DQ17, HLZ20} and switching graphs \cite{CD19, DHJK20, LX10}, variants of the CS model \cite{DHK19, HKR20, MT11}, etc. We also refer to survey papers \cite{CHZ17, MM19} for more references.

In this paper, we are interested in the interplay between singular communication weights 
\[
\psi(r) := r^{-\alpha}, \quad \alpha>0.
\]
and general network topologies. Under the all-to-all interaction network, i.e. $a_{ij}=1$ for all $i,j = 1, \dots, N$, a notable feature of singular communication weights is that when $\alpha \ge1$, agents governed by \eqref{main_sys} do not collide for all time if they are initially non-collisional \cite{CCMP17}, which ensures the global-in-time existence of smooth solutions to \eqref{main_sys}. However, when $\alpha \in (0,1)$, collisions between agents may occur in a finite time, leading to the breakdown of smooth solutions \cite{P14}. 

For regular communication weights, on the other hand, a sufficient framework for the emergence of flocking is establised when there is a spanning tree in the network topology \cite{DQ17}. In particular, in one-dimensional case, if the communcation weight is non-integrable, the existence of a spanning tree guarantees the emergence of mono-cluster flocking for any initial configuration \cite{HLZ20}.

Motivated by these observations, we address the following two questions:

\begin{itemize}
\item
(Global-in-time existence problem) Under what conditions on the network topology do smooth solutions to \eqref{main_sys} exist globally in time? In particular, does the condition $\alpha\geq 1$ alone guarantee global-in-time existence, as in the all-to-all case? If not, which network structures ensure global-in-time existence?

\item
(Emergence of mono-cluster flocking)  Assuming the global-in-time existence of smooth solutions, what conditions on the network topology and initial configuration guarantee the emergence of mono-cluster flocking in \eqref{main_sys}?
\end{itemize}

Our main results are two-fold. First, we establish the global-in-time existence of smooth solutions to \eqref{main_sys} for several classes of network topologies when $\alpha \ge 1$. Here, we remark that the global existence  is not necessarily equivalent to the absence of collisions between agents. Indeed, if two agents $i$ and $j$ collide at some time $t=t_0$ without interacting with each other, i.e. $a_{ij} = a_{ji} = 0$, their collision does not induce a singularity in the equations of $\dot{v}_{i}$ and $\dot{v}_j$, unless they collide with other agents. Thus, in terms of the global existence, we only focus on collisions between interacting agents and distinguish them from those between noninteracting agents (see Definition \ref{def:col}). 

In one-dimensional setting, we further characterize the network topologies under which system \eqref{main_sys} admits a global-in-time smooth solution for every non-collisional initial configuration. Conversely, for network topologies not belonging to this class, we can construct a non-collisional initial data whose corresponding solutions experience finite-time collisions. Especially, it shows that the condition $\alpha\geq 1$ alone does not guarantee collision avoidance under general network topologies, in contrast to the all-to-all case.

Second, we present the sufficient conditions for the emergence of mono-cluster flocking in system \eqref{main_sys}. As in the case of regular communication weights, we establish a sufficient framework for the emergence of mono-cluster flocking when the network topology has a spanning tree. We further investigate the case where the network topology is scrambling. Under one of the assumptions for the global-in-time existence theory, we show that when $\alpha=1$, system \eqref{main_sys} exhibits mono-cluster flocking for any non-collisional initial configuration.

The rest of this paper is as follows. In Section \ref{sec:pre}, we provide concepts on directed graphs and technical estimates used in the paper. In Section \ref{sec:main_result}, we present our main results on the global well-posedness and mono-cluster flocking estimates. In Section \ref{sec:glob_ex}, we establish the global-in-time existence theory of smooth solutions to system \eqref{main_sys}. In Section \ref{sec:flock}, we discuss the sufficient framework for the emergence of mono-cluster flocking. Finally, Appendix \ref{sec:app.A}
is devoted to the characterization of one of the class of digraphs described in Theorem \ref{thm:glob_ex}.

\vspace{0.2cm}

\noindent {\bf Notation.} Throughout this paper, we write $X= (x_1,...,x_N) \in (\R^d)^N$ and $V = (v_1,...,v_N) \in (\R^d)^N$.  $\Gamma = \Gamma(t)$ and $\Lambda=\Lambda(t)$ stand for the diameters for position and velocity:
\[\Gamma(t) = \displaystyle\max_{1\le i,j \le N} |x_i(t)-x_j(t)|, \quad \mbox{and} \quad \Lambda(t) = \displaystyle\max_{1\le i,j \le N} |v_i(t)-v_j(t)|.\]

For $m\times n$ real matrices $A = (a_{ij})$ and $B = (b_{ij})$, we write $A \ge B$ if $a_{ij} \ge b_{ij}$ for every $i$ and $j$. We also denote by $\mathbb{I}_n$ the $n \times n$ identity matrix for each $n \in \N$.

\vspace{0.2cm}

%
%
%
%

\section{Preliminaries}\label{sec:pre}
\setcounter{equation}{0}

In this section, we  study basic properties of solutions to \eqref{main_sys} and recall several concepts on directed graphs, scrambling matrices and state-transition matrices.

\subsection{Technical estimates}
Here, we present several technical estimates used throughout this paper. We first provide basic estimates for the velocity and position diameters.

\begin{proposition}\label{prop:dm_bd}
    For $T>0$, let $(X,V)$ be a solution to \eqref{main_sys} on the time interval $[0,T)$. Then we have
    \[\Lambda(t) \le \Lambda(0) \quad \mbox{and} \quad \Gamma(t) \le  \Gamma(0) + t\Lambda(0)\]
    for every  $t \in [0,T]$.
    \end{proposition}
\begin{proof} 

For $t\in(0,T)$ at which $\Lambda$ is differentiable, let $\Lambda(t) = |v_i(t)-v_j(t)|$.
\[
    \frac{1}{2}\frac{d}{dt}\Lambda^2 = \frac{1}{N}\sum_{k\neq i} a_{ik} \psi(|x_k-x_i|)(v_k-v_i)\cdot(v_i-v_j) - \frac{1}{N}\sum_{k\neq j} a_{jk} \psi(|x_k-x_j|)(v_k-v_j)\cdot(v_i-v_j)
\]
By the Cauchy-Schwarz inequality and definition of $\Lambda$, we derive that
\begin{eqnarray*}
    (v_k-v_i)\cdot(v_i-v_j) &=& (v_k-v_j)\cdot(v_i-v_j) - |v_i-v_j|^2 \\
    &\le& |v_i-v_j| (|v_k-v_j|-|v_i-v_j|) \le 0 \\\\
    -(v_k-v_j)\cdot(v_i-v_j)&=& -|v_i-v_j|^2 + (v_i-v_k)\cdot(v_i-v_j)\\
    &\le&- |v_i-v_j|(|v_i-v_j| - |v_i-v_k|) \le 0
\end{eqnarray*}
Since $\psi \ge 0$, $\displaystyle\frac{d}{dt} \Lambda\le 0$ a.e. $t\in (0,T)$.
Thus, $\Lambda(t) \le \Lambda(0)$ for every $t \in [0,T]$.

\noindent Similarly, we obtain $\displaystyle \frac{d}{dt}\Gamma(t)\le \Lambda(t) \le \Lambda(0)$ a.e. $t \in (0,T)$. Therefore, $\Gamma(t) \le \Gamma(0) + \Lambda(0) t$ for every $t\in[0,T]$. 

\end{proof}

Next, we provide a technical lemma which will be used for the existence proof.

\begin{lemma}\label{lem:ineq_ex}
    Let $t_0>0$, $f: (0,\infty) \to (0,\infty)$ be continuous, and $\phi$ be a positive locally Lipschitz continuous function on the interval $[0,t_0)$. Suppose that there exists a nonnegative locally Lipschitz continuous and bounded function $\varphi$ on $[0,t_0)$ such that for a.e. $t\in(0,t_0)$, 
    \[
    \left|\dot\phi(t)\right| \le c_0 \varphi(t) \quad \mbox{and} \quad \dot\varphi(t)\le  -c_1(f\circ\phi)(t)\varphi(t)+c_2
    \]
     for some constants $c_0,c_1,c_2>0$. Then $F \circ \phi$ is bounded on $[0,t_0)$, where $F$ is a primitive of $f$. If additionally $\displaystyle \lim_{r\to 0^+}F(r) = -\infty$, then $\displaystyle\inf_{t\in [0,t_0)} \phi(t) >0$ and $\displaystyle \lim_{t\to t_0^-}\phi(t) =: \phi(t_0) > 0$.
\end{lemma}
\begin{proof} 
Although similar proofs can be found in \cite{CCMP17}, we present it for readers' convenience. Since $\varphi$ is bounded, there exists $M>0$ such that, for $t \in [0,t_0)$,
\begin{eqnarray*}
    \int_0^tc_0(f\circ \phi)(\tau)\varphi(\tau)d\tau&\le& \frac{c_0}{c_1}\int_0^t c_2-\dot\varphi(\tau)d\tau\\
    &\le&\frac{c_0}{c_1}  (c_2t_0-\varphi(t) + \varphi(0))\\
    &\le& M.
\end{eqnarray*}
For $t \in (0,t_0)$ at which $F\circ \phi$ is differentiable, 
\[\left|\left(\frac{d}{dt}(F\circ \phi)\right)(t)\right|\le (f\circ\phi)(t)|\dot \phi(t)|\le c_0 (f\circ\phi)(t) \varphi(t),\]
so, $F\circ \phi$ is bounded on $[0,t_0)$. Assume that $\displaystyle\lim_{r \to 0^+} F(r) = - \infty$ and $\displaystyle \inf_{t \in [0 ,t_0)} \phi(t) =0$. 

Then $\displaystyle \inf_{t \in [0,t_0)} (F\circ \phi)(t) = -\infty$. Thus, $\displaystyle\inf_{t \in [0,t_0)}\phi(t) >0$. Since $\phi$ is globally Lipschitz by $|\dot \phi| \le c_0 \varphi$, $\displaystyle \phi(t_0) > 0$. \hfill

\end{proof}

\subsection{Directed graphs}
In this subsection, we review several concepts on directed graphs (digraphs). A digraph $G =(V_G,E_G)$ consists of a set of vertices (or nodes) $V_G = \{1, \dots, N\}$, and a set of edges (or arcs) $E_G \subseteq V_G \times V_G$ comprising ordered pairs of vertices. Here, we say a vertex $j$ is a {\bf neighbor} or {\bf predecessor} of vertex $i$ if $(j,i) \in E_G$ and we also write it as $j \to i$. If $(i, i) \in E_G$, we say $G$ has a {\bf self-loop} at vertex $i$. For a given digraph $G= (V_G, E_G)$, we consider its 0-1 adjacency matrix $A=(a_{ij})$:
\[
a_{ij} = \left\{\begin{array}{lcl} 1 & \mbox{if} & (j, i) \in E_G,\\
0 &\mbox{if} & (j,i) \notin E_G.\end{array} \right.
\]

A {\bf path} in $G$ from $j$ to $i$ is a sequence of ordered, distinct vertices $j_0 = j$, $j_1$, $\dots$, $j_n =i$ such that $(j_{m-1}, j_m) \in E_G$ for every $1 \le m \le n$, or equivalently,
\[
j_0 = j \to j_1 \to \cdots \to j_n = i.
\]
If there is a path from $j$ to $i$, we say $i$ is {\bf reachable} from $j$. In this case, we define the distance $d(j,i)$ from $j$ to $i$ as the length of the shortest path from $j$ to $i$.  Furthermore, we say a digraph $G$ has a {\bf spanning tree} if $G$ has a vertex $j$ from which any other vertices are reachable. Such a vertex is called a {\bf root} of $G$.

Below, we present several concepts on digraphs which will be crucially used throughout the paper:

\begin{definition}\label{def:dig}
Let $G = (V_G,E_G)$ be a digraph and $A = (a_{ij})$ be its 0-1 adjacency matrix. We say \\

\begin{enumerate}
\item $G$ is {\bf semicomplete} if

\begin{center}
     $\forall i,j \in V_G$, $a_{ij} + a_{ji}>0,$
\end{center}

\item
 $G$ is {\bf scrambling} if
\begin{center}
    $\forall i,j \in V_G$, $\exists r=r(i,j) \in V_G$ such that $a_{ir}a_{jr} = 1,$
\end{center}

\item
$G$ is an {\bf out-arborescence} if $G$ has a unique vertex $r$ such that, after removing all the self-loops in $G$, there is exactly one path from $r$ to $i$  for every $i \in V_G\setminus \{r\}$,

\vspace{0.2cm}

\item
 $G$ is an {\bf undirected graph} if 
\begin{center}
     $\forall i,j \in V_G$, $a_{ij} = a_{ji}$.
\end{center}

\end{enumerate}

\end{definition}

\begin{remark}\label{rm:scrb}
Let $B = (b_{ij})$ be a nonnegative $N\times N$ matrix. Its ergodicity coefficient is defined as
\[
\mu(B) := \min_{i,j} \sum_{k=1}^N \min\{b_{ik}, b_{jk}\}.
\]
Then we have the followings:
\begin{enumerate}
\item
If $B_1$ and $B_2$ are nonnegative matrices, then
\[
B_1 \ge B_2 \quad \Longrightarrow \mu(B_1) \ge \mu(B_2).
\]
\item
Suppose that $B$ is an 0-1 adjacency matrix of a digraph $G = (V_G, E_G)$ which has a spanning tree and a self-loop at each vertex. Then one has
\[
\mu(B^h) \ge 1,
\]
where $h$ is the shortest height of all spanning trees of $G$, i.e.
\[
h:= \min \lt\{\ \max_{i \in V_G} d(r,i) \ : \ r \mbox{ is a root of } G \rt\}.
\]
See \cite[Proposition 1]{DQ17} for its proof. 
\end{enumerate}

\end{remark}

We also present the concept of induced subgraphs and connected digraphs, which will appear in the statement of our main results.

\begin{definition}\label{def:subgr}
    Let $G =(V_G,E_G)$ be a digraph and $A$ be a subset of $V_G$. The \textbf{induced subgraph} $G[A]$ is the digraph with the vertex set $A$ and the edge set consisting of the edges in $E_G$ that have both tail and head in $A$, i.e.
    \begin{center}
        $G[A] = (A, E[A])$ where $E[A] = \{(i,j)\in E_G : i,j \in A\}$.
    \end{center}
\end{definition}

\begin{remark}\label{rmk:sc_sc}

Let $G$ be a digraph and $A$ be a subset of $V$. We can easily deduce the followings:
\begin{enumerate}
\item
If $G$ is semicomplete, then $G[A]$ is also semicomplete.

\item
Suppose that each vertex of $G$ has a self-loop. Then $G$ is scrambling if $G$ is semicomplete. Furthermore, the above implies that  every induced subgraph of $G$ is scrambling if $G$ is semicomplete.
\item
Conversely, one can prove that if every induced subgraph of $G$ is scrambling, then $G$ is semicomplete.
\end{enumerate}
\end{remark}
Finally, we provide the concept of connectivity.

\begin{definition}\label{def:cnt_gp}
    Let $G =(V_G,E_G)$ be a digraph. The \textbf{underlying graph} $G^*=(V_G,E_G^*)$ of $G$ is the graph such that $E_G^* = \{(i,j) \mid (i,j) \in E_G\text{ or }(j,i) \in E_G\}$. If the underlying graph $G^*$ has a spanning tree, $G$ is called \textbf{weakly connected}. A \textbf{weakly connected component} of $G$ is a maximal weakly connected induced subgraph. In other words, a subset $A$ of $V_G$ induces a weakly connected component of $G$ if $G[A]$ is weakly connected and there is no subset $B$ of $V_G$ such that $A\subsetneq B$ and $G[B]$ is weakly connected. 
\end{definition}
\begin{remark}
Unlike \textbf{strong connectivity}, which requires a directed path between every ordered pair of vertices, weak connectivity only requires the underlying graph to be connected. 
\end{remark}

\subsection{State-transition matrix}

In this subsection, we recall the concept of state-transition matrices and recast system \eqref{main_sys} in the state-transition matrix form. Let $t_0 \in \R$ and $B:[t_0, \infty) \to \R^{N\times N}$ be an $N\times N$ matrix of smooth functions. Then we consider the initial-value problem associated to $B$:
\bq\label{eq:sys_ode}
\begin{aligned}
\frac{d}{dt}\Xi(t) &= B(t)\Xi(t), \quad t >t_0,\\
\Xi \Big|_{t=t_0} &= \Xi(t_0).
\end{aligned}
\eq
Then the unique solution to \eqref{eq:sys_ode} is given by
\[
\Xi(t) = \Phi(t, t_0) \Xi(t_0), \quad t\ge t_0,
\]
where $\Phi = \Phi(t,t_0)$ is called the state-transition matrix for \eqref{eq:sys_ode} and it can be written in the form of Peano-Baker series (see \cite{S98}):
\[
\Phi(t,t_0) = \mathbb{I}_N + \sum_{k=1}^\infty\int_{t_0}^t \int_{t_0}^{\tau_1}\cdots\int_{t_0}^{\tau_{n-1}}A(\tau_1)A(\tau_2)\cdots A(\tau_n)\,d\tau_n \cdots d\tau_2 d\tau_1.
\]
We also consider two initial-value problems
\[
\frac{d}{dt}\Xi(t) = B(t)\Xi(t), \quad \frac{d}{dt}\Xi(t) = (B(t) + c\mathbb{I}_N) \Xi(t), \quad t > t_0,
\]
and let $\Phi(t,t_0)$ and $\Psi(t,t_0)$ be the state-transition matrices corresponding to the above systems, respectively. Then we can find a relation between these state-transition matrices:
\bq\label{eq:stat_tran}
\Phi(t, t_0) = e^{-c(t-t_0)}\Psi(t,t_0), \quad t \ge t_0
\eq
and we refer to \cite{DHK19} for its proof.

\vspace{0.2cm}

%
%
%
%

\section{Main results}\label{sec:main_result}
\setcounter{equation}{0}
In this section, we present our main results. Throughout the paper, we impose the following standing assumptions on the digraph representing the network topology in \eqref{main_sys}:

\begin{itemize}
\item[$(\mathcal{A}1)$] $G$ is weakly connected.
\item[$(\mathcal{A}2)$] Every vertex in the digraph $G= (V_G, E_G)$ representing the network topology has a self-loop.
\end{itemize}

\vspace{0.2cm}

These assumptions do not involve any loss of generality. Indeed, if $G$ is not weakly connected, we can decompose system \eqref{main_sys} into independent subsystems corresponding to the weakly connected components of $G$, and we can apply our argument to each subsystem separately. Moreover, since any diagonal entries of 0-1 adjacency matrix $A = (a_{ij})$ do not appear in \eqref{main_sys}, presence of any self-loops do not change the dynamics.

As we emphasized in the introduction, it is desirable to distinguish collisions between interacting agents from those between noninteracting agents. This motivates the following definition.

\begin{definition}\label{def:col}
Let $(X,V)$ be a solution to \eqref{main_sys} and $A=(a_{ij})$ be the associated 0-1 adjacency matrix. We say agents $i$ and $j$ \textbf{overlap} at $t=t_0$ if $x_i(t_0) = x_j(t_0)$. Moreover, we say they \textbf{collide} at $t=t_0$ if they overlap at $t=t_0$ and $a_{ij} +a_{ji} >0$. 
\end{definition}

Then our first result on the global existence is as follows:

\begin{theorem}\label{thm:glob_ex}
Let $\alpha \ge 1$ and  $G= (V_G, E_G)$ be the digraph representing the network topology in \eqref{main_sys} with 0-1 adjacency matrix $A= (a_{ij})$. Suppose that the initial data $(X(0), V(0))$ are non-collisional, i.e.
\[x_i(0) \neq x_j(0) \quad \mbox{for} \quad 1\le i \neq j \le N \quad \mbox{with} \quad a_{ij} + a_{ji} >0.\]
Assume that one of the followings hold:

\begin{enumerate}
\item
Every weakly connected, induced subgraph of $G$ is scrambling,
\item
Every weakly connected, induced subgraph of $G$ is an out-arborescence,
\item
$G$ is undirected.
\end{enumerate}
Then, system \eqref{main_sys} admits a global-in-time smooth solution. Moreover, agents governed by \eqref{main_sys} possibly overlap, but do not collide in a finite time, i.e.
\[
x_i(t) \neq x_j(t) \quad \mbox{for}\quad 1\le i \neq j \le N \quad \mbox{with} \quad a_{ij} + a_{ji} >0,
\]
for any $t \ge 0$.

\end{theorem}


We next turn to the one-dimensional case.  Before stating our result, we introduce the following concept.

\begin{definition}
    Let \(G=(V_G,E_G)\) be a digraph. A vertex \(s\in V_G\) with in-degree zero is called a \textbf{source} of \(G\), i.e. 
\[\
(i,s)\notin E_G \quad \mbox{for all} \quad  i\in V_G\setminus\{s\},
\]
or equivalently, 
\[
a_{si} =0 \quad\mbox{for all}\quad i \in V_G \setminus \{s\}.
\]
\end{definition}
We are now ready to present the one-dimensional case.

\begin{theorem}\label{thm:1d_glob}
    Let $d=1$ and $G= (V_G, E_G)$ be the digraph representing the network topology in \eqref{main_sys} with 0-1 adjacency matrix $A= (a_{ij})$. Then we have the followings:
        \begin{itemize}
        \item[(i)] Assume that $\alpha \ge 1$ and every connected induced subgraph of $G$ has at most one source. Then for every non-collisional initial configuration $(X(0), V(0))$, i.e.
\[x_i(0) \neq x_j(0) \quad \mbox{for} \quad 1\le i \neq j \le N \quad \mbox{with} \quad a_{ij} + a_{ji} >0,\]

    system \eqref{main_sys} admits a global-in-time smooth solution. Moreover, agents governed by \eqref{main_sys} possibly overlap, but do not collide in a finite-time, i.e.
\[
x_i(t) \neq x_j(t) \quad \mbox{for}\quad 1\le i \neq j \le N \quad \mbox{with} \quad a_{ij} + a_{ji} >0,
\]
for any $t \ge 0$.
\vspace{0.1cm}

\item[(ii)] Assume that $\alpha >0$ and there exists a connected induced subgraph of $G$ has at least two sources. Then we can construct an initial data $(X(0), V(0))$ such that the corresponding smooth solution $(X,V)$ to \eqref{main_sys} can exist up to a finite time, i.e. the solution experiences a finite-in-time collision.
    \end{itemize}
\end{theorem}

\vspace{0.2cm}

\begin{remark} We leave some comments on the results of Theorems \ref{thm:glob_ex} and \ref{thm:1d_glob}.

\begin{enumerate}
\item
We deduce from Remark \ref{rmk:sc_sc} and Theorem \ref{thm:glob_ex} that if $G$ is semicomplete and the initial data are non-overlapping, i.e.
\[x_i(0) \neq x_j(0) \quad \mbox{for} \quad 1\le i \neq j \le N,\]
then agents governed by \eqref{main_sys} do not overlap in a finite time.
\item
One directly deduces from the result that the sole condition $\alpha \ge 1$ does not guarantee the global-in-time existence of smooth solutions to \eqref{main_sys}. Note that in \cite{CCMP17}, only the all-to-all network was considered and in this case there is no source in the network topology.

\item
Note that the condition in Theorem \ref{thm:1d_glob} can be fulfilled by those in Theorem \ref{thm:glob_ex}. Clearly, undirected graphs and out-arborescences have at most one source. Moreover, the same holds for scrambling graphs. Indeed, suppose there exist two distinct sources $s_1$ and $s_2$. By definition of the source, there is no vertex $i$ satisfying both $(i,s_1) \in E$ and $(i,s_2)\in E$. Hence, the global-in-time existence results in the two theorems are consistent. 
\end{enumerate}
\end{remark}

Next, we proceed to our main results on the emergent behavior of \eqref{main_sys}. We hereby recall the definition for the mono-cluster flocking.

\begin{definition}\label{def:flock}
    For $X= (x_1,...,x_N)$ and $V = (v_1,...,v_N)$, let 
    \[\Lambda(t) = \displaystyle\max_{1\le i,j \le N} |v_i(t)-v_j(t)| \text{\quad and\quad} \Gamma(t) = \displaystyle\max_{1\le i,j \le N} |x_i(t)-x_j(t)|.\]
    We say that system \eqref{main_sys} exhibits \textbf{mono-cluster flocking} if the following conditions are satisfied:
    \[\lim_{t \to\infty} \Lambda(t) =0 \text{\quad and \quad} \sup_{t\ge0} \Gamma(t) < \infty.\]
\end{definition}

Let $G= (V_G, E_G)$ be the digraph representing the network topology in \eqref{main_sys}, we define the quantity $\gamma = \gamma(t)$ as
\[
\gamma(t) := \min\{ \ |(x_i - x_j)(t)| \ : \ (i,j) \in E_G, \quad i \neq j\}.
\]
Then the global-in-time existence of smooth solutions to \eqref{main_sys} is equivalent to $\gamma(t) >0$ for every $t \ge 0$. Now, we are ready to provide our first result on the mono-cluster flocking.

\begin{theorem}\label{thm:LT1}
Let $\alpha>0$ and $G= (V_G, E_G)$ be the digraph representing the network topology in \eqref{main_sys} with a spanning tree. Suppose that the initial data $(X(0), V(0))$ are non-collisional, i.e. $\gamma(0)>0$, and there exist $\rho >0$ and $\eta \in (0,\gamma(0))$ such that
    \[\Lambda(0) \le \frac{\psi_\eta \delta}{h} \min\{\rho,\eta\},\]
    where $\delta = \displaystyle\frac{1}{h!}\left(\frac{h\psi_\rho}{Ne\psi_\eta}\right)^h$, $\psi_\rho = \psi(\Gamma(0)+\rho)$, $\psi_\eta = \psi(\gamma(0)-\eta)$, and $h$ is the shortest height of all spanning trees of $G$. Then system \eqref{main_sys} admits a global-in-time smooth solution and exhibits mono-cluster flocking.
\end{theorem}

Combined with global-in-time existence theory developed in Theorem \ref{thm:glob_ex}, we also derive the following result:
\begin{theorem}\label{thm:LT2}
Let $\alpha>0$ and $G= (V_G, E_G)$ be the digraph representing the network topology in \eqref{main_sys}. Suppose that the initial data $(X(0), V(0))$ are non-collisional,  i.e. $\gamma(0)>0$, and one of the following conditions holds:
    \begin{enumerate}
        \item $\alpha \in (0,1)$, $G$ is scrambling, and 
        \[\Lambda(0)<\frac{S(G)\gamma(0)}{N \Gamma_*^\alpha}\]
        \item $\alpha>1$, every connected induced subgraph of $G$ is scrambling, and
        \[\Lambda(0)<\frac{S(G)\Gamma(0)^{1-\alpha}}{N(\alpha-1)}, \]
        \item $\alpha=1$ and every connected induced subgraph of $G$ is scrambling,
    \end{enumerate}
    where $S(G)= \displaystyle \min_{i\neq j} |\{k \in V_G : k \to i\}\cap \{k \in V_G : k \to j\}|$ and $\Gamma_*$ is defined as
    \[\Gamma_* = \lt\{\begin{array}{lcl}
        \displaystyle\Gamma(0)\exp\left(\frac{N\Lambda(0)}{S(G)}\right)& \mbox{if} & \alpha =1,\\
        \displaystyle \left(\Gamma(0)^{1-\alpha}+\frac{(1-\alpha)N} {S(G)}\Lambda(0)\right)^{1/(1-\alpha)}& \mbox{if} & \alpha \in (0,1) \cup (1, \infty).
    \end{array}\rt.\]
     Then system \eqref{main_sys} admits a global-in-time smooth solution and exhibits mono-cluster flocking:
    \[\Gamma(t) \le \Gamma_*, \quad \Lambda(t) \le  \Lambda(0)\exp\left(-\frac{S(G)}{N\Gamma_*^\alpha}t\right) \quad \mbox{for all} \quad t\ge 0.\]

\end{theorem}

\begin{remark}
The estimates in Theorem \ref{thm:LT1} require the communication weight between interacting agents to be uniformly bounded in time. Thus, one needs a uniform-in-time lower bound on the  distances between interacting agents, which is not guaranteed by Theorem \ref{thm:glob_ex}. Accordingly, conditions in Theorem \ref{thm:LT1} involves the initial minimal distance for the emergence of flocking. In contrast, the estimates in Theorem \ref{thm:LT2} do not require the boundedness of communication weights, but only a priori existence of smooth solutions. However, since no network topology guarantees global-in-time existence for $\alpha\in(0,1)$, the initial minimum distance condition appears only in (i) of Theorem \ref{thm:LT2}. In particular, under condition (i), we can obtain
\[
\inf_{t \ge 0}\gamma(t) >0.
\]

\end{remark}

%
%
%
%

\section{Global-in-time existence of smooth solutions}\label{sec:glob_ex}
\setcounter{equation}{0}

In this section, we provide the proof of Theorems \ref{thm:glob_ex} and \ref{thm:1d_glob}. We start with Theorem \ref{thm:glob_ex}, which covers several classes of network topologies. Accordingly, we consider each case separately below:

\subsection{Proof of Theorem \ref{thm:glob_ex} under condition (1)}\label{ssec:glob_ex_pf1}
We establish the global existence of smooth solutions to \eqref{main_sys} when every connected, induced subgraph of $G$ is scrambling,

\vspace{0.2cm}

Under our assumption on the initial configuration, the local-in-time existence of a unique smooth solution is straightforward. Thus, there exists a unique smooth solution exists on some time interval $[0,T)$. Then we have two possibilities: either any of the agents do not collide up to $t=T$, or a collision, not overlapping, occurs at some time $t_0 \in(0,T]$. Assume for a contradiction that such $t_0$ exist. Then we may assume that $t=t_0$ is the first collision time and the solution uniquely exists and remains smooth on $[0,t_0)$.

By definition, there exists an index $\ell\in V_G$ such that the $\ell$-th agent collides with another agent. We collect the indices of agents which collide with the $\ell$-th agent  to form $C_\ell \subseteq V_G$ and if necessary, we expand $C_\ell$ further so that $G[C_\ell]$ is a weakly connected component of $G$.  By our assumption, $G[C_\ell]$ is scrambling. Then, we can find $\delta>0$ satisfying
\[|x_i(t)-x_j(t)|\ge \delta > 0\; \text{ in}\; [0, t_0) \quad \mbox{for all} \quad i \in C_\ell, \quad j \notin C_\ell \quad \mbox{with} \quad a_{ij}=1. \]

For simplicity, we write 
\[
x_{ij} := x_i-x_j, \quad v_{ij} := v_i-v_j,\quad\mbox{and}\quad  \psi_{ij} := \psi(|x_{ij}|).
\] We also define $D_X[C_\ell]$ and $ D_V[C_\ell]$ as follows:
\begin{center}
    $\displaystyle D_X[C_\ell](t) := \sqrt{\max_{i,j \in C_\ell} |x_{ij}(t)|^2}$ \;\;and\;\; $\displaystyle D_V[C_\ell](t) := \sqrt{\max_{i,j \in C_\ell} |v_{ij}(t)|^2},$
\end{center}
and our assumption implies $\displaystyle D_X[C_\ell](t) \to 0$ as $t \nearrow t_0$.

Then we claim that the following differential inequalities hold for a.e. $t\in (0,t_0)$:
    \[
        \left|\frac{d}{dt}  D_X[C_\ell]\right| \le  D_V[C_\ell], \quad \mbox{and}
    \]
\[
      \frac{d}{dt}  D_V[C_\ell] \le -c_0\psi( D_X[C_\ell]) D_V[C_\ell] + c_1\quad \mbox{for some} \quad  c_0,c_1>0.
\]
Once we establish the above inequalities, we can use Lemma \ref{lem:ineq_ex} to derive a contradiction.\\ 

\noindent Now, for $t\in(0,t_0)$ at which $ D_X[C_\ell]$ and $ D_V[C_\ell]$ are differentiable, we can easily deduce
\[   
 \left|\frac{d}{dt}  D_X[C_\ell]\right| \le  D_V[C_\ell].
\]

At time $t$, we may let $ D_V[C_\ell] =|v_{ij}|$ for some $i,j \in C_\ell$. Then

\[
\begin{aligned}
    \frac12\frac{d}{dt}  D_V[C_\ell]^2 &= \frac{1}{N} \sum_{k = 1}^N\left(a_{ik}\psi_{ik}v_{ki}\cdot v_{ij}-a_{jk}\psi_{jk}v_{kj}\cdot v_{ij}\right)\\
    &=\frac{1}{N}\left(\sum_{k \in C_\ell} + \sum_{k\not \in C_\ell }\right)\left(a_{ik}\psi_{ik}v_{ki}\cdot v_{ij}-a_{jk}\psi_{jk}v_{kj}\cdot v_{ij}\right)\\
    &=: \mathcal{J}_1+\mathcal{J}_2.
\end{aligned}
\]
For $\mathcal{J}_1$, since the maximality of $i$ and $j$ implies
\[
 v_{kj}\cdot v_{ij} \ge 0, \quad \forall k \in C_\ell,
\]

we get

\[
\begin{aligned}
    \mathcal{J}_1 &= -\frac{1}{N} \sum_{k \in C_\ell}\left(a_{ik}\psi_{ik}v_{ik}\cdot v_{ij}+a_{jk}\psi_{jk}v_{kj}\cdot v_{ij}\right)\\
    &\le -\frac{1}{N} \psi( D_X[C_\ell])\sum_{k \in C_\ell}\left(a_{ik}v_{ik}\cdot v_{ij}+a_{jk}v_{kj}\cdot v_{ij}\right)
\end{aligned}
\]
Since $G[C_\ell]$ is scrambling, there exists $r \in C_\ell$ such that $a_{ir} = a_{jr} =1$. This implies
\[
    \sum_{k \in C_\ell}\left(a_{ik}v_{ik}\cdot v_{ij}+a_{jk}v_{kj}\cdot v_{ij}\right) \ge a_{ir} v_{ir} \cdot v_{ij} + a_{jr}v_{rj}\cdot v_{ij} = |v_{ij}|^2. 
    \]
Hence, we can obtain

\[
    \mathcal{J}_1 \le -\frac{1}{N} \psi( D_X[C_\ell]) |v_{ij}|^2 = -\frac{1}{N} \psi( D_X[C_\ell])  (D_V[C_\ell])^2.
\]
For $\mathcal{J}_2$, we use $|x_{ik}|, |x_{jk}| \ge \delta$ for $k \notin C_\ell$ and Proposition \ref{prop:dm_bd} to yield

\[
\begin{aligned}
    \mathcal{J}_2 &\le  \frac{\psi(\delta)}{N}  \sum_{ k\not \in C_\ell} ( |v_{ki}| + |v_{kj}|) |v_{ij}|\\
    &\le 2\psi(\delta) \Lambda(0) D_V[C_\ell]
\end{aligned}
\]

\noindent Hence, we collect the estimates for $\mathcal{J}_1$ and $\mathcal{J}_2$ to obtain
\[
    \frac{d}{dt} D_V[C_\ell] \le -\frac{1}{N}\psi( D_X[C_\ell]) D_V[C_\ell] +2\psi(\delta) \Lambda(0),
  \]
which proves the claim. Since the primitive $\Psi$ of $\psi$ is given as 
\bq\label{eq:prim_psi}
\Psi(r) =
\begin{cases}
        \log r&\text{if } \alpha =1,\\
        \displaystyle\frac{1}{1-\alpha} r^{1-\alpha}&\text{if } \alpha >1.
\end{cases}
\eq
we have $\displaystyle\lim_{r \to 0^+}\Psi(r) = - \infty$. Hence, we deduce from Lemma \ref{lem:ineq_ex}  that 
\[
\lim_{t \to t^-_0}D_X[C_\ell](t) >0,
\]
which contradicts to $D_X[C_\ell](t) \to 0$ as $t \to t_0$. This completes the proof of Theorem \ref{thm:glob_ex} under condition (1).

\begin{remark}
The class of digraphs satisfying the condition (1) in Theorem \ref{thm:glob_ex} might seem unfamilar. Interestingly, however, each digraph in this class can be generated recursively from semicomplete digraphs by applying certain operations to pairs of digraphs. Although this provides an interesting characterization, it lies somewhat outside the scope of the present work. For interested readers, we present the details in Appendix \ref{sec:app.A}.
\end{remark}

\vspace{0.2cm}

\subsection{Proof of Theorem \ref{thm:glob_ex} under condition (2)}
Next, we assume that every connected induced subgraph of $G$ is out-arborescence.

As in the proof of Theorem \ref{thm:glob_ex} under condition (1), we assume for a contradiction that a collision occurs in a finite time. Then we set $t_0$, $C_\ell$, and $\delta>0$ as before.

Now, we define the following functional to measure the distance between agents in $C_\ell$:
\[\mathcal{P}_{C_\ell}(t) = \sum_{\substack{i,j\in C_\ell\\ i\neq j}}a_{ij}\Psi(|x_{ij}(t)|)\]
where $\Psi$ is the primitive of $\psi$ defined in \eqref{eq:prim_psi}.

Then $\mathcal{P}_{C_\ell}(t) \to -\infty$ as $t \nearrow t_0$. Since $G[C_\ell]$ is weakly connected, it is an out-arborescence. So $G [C_\ell]$ has a unique vertex $r$ such that, after removing all the self-loops in $G[C_\ell]$, there is exactly one path from $r$ to $i$  for every $i \in C_\ell\setminus \{r\}$.
This implies that for every $i \in C_\ell\setminus\{r\}$, there exists a unique vertex $ j\in C_\ell\setminus\{i\}$ such that $a_{ij}=1$. Denote such $j$ of $i$ as $\pi(i)$. Note that $a_{\pi(i),i}=0$ for all $i \in C_\ell\setminus\{r\}$. Then we can rewrite $\mathcal{P}_{C_\ell}$ as

\[\mathcal{P}_{C_\ell}(t) = \sum_{i \in C_\ell\setminus\{r\}}\Psi(|x_{i,\pi(i)}(t)|).\]
One easily deduces that
\[\left|\frac{d}{dt}\mathcal{P}_{C_\ell}\right|= \left|\sum_{i\in C_{l}\setminus\{r\}} \psi_{i,\pi(i)}\frac{x_{i,\pi(i)}\cdot v_{i,\pi(i)}}{|x_{i,\pi(i)}|}\right|\le \sum_{i\in C_{l}\setminus\{r\}}\psi_{i,\pi(i)} |v_{i,\pi(i)}|.\]
Thus,

\[|\mathcal{P}_{C_\ell}(t)|\le |\mathcal{P}_{C_\ell}(0)| + \sum_{i\in C_\ell\setminus\{r\}}\int_0^t\psi_{i,\pi(i)}(\tau)|v_{i,\pi(i)}(\tau)|d\tau\text{ for all }t\in [0,t_0).\]
We will proceed by induction on the distance to the unique root $r$ to show $\displaystyle\int_0^t\psi_{i,\pi(i)}(\tau)|v_{i,\pi(i)}(\tau)|d\tau$ is bounded on $[0,t_0)$ for every $i \in C_\ell\setminus\{r\}$.\\

\vspace{0.1cm}

\noindent $\bullet$  (Initial step) Let $i\in C_\ell\setminus \{r\}$ such that $\pi(i)=r$, i.e. $d(r,i)=1$. Then we have
\[
\begin{aligned}
    \frac12\frac{d}{dt} |v_{ir}|^2 &=\frac{1}{N}\left(\sum_{k \in C_\ell} + \sum_{k\not \in C_\ell }\right)\left(a_{ik}\psi_{ik}v_{ki}\cdot v_{ir}-a_{rk}\psi_{rk}v_{kr}\cdot v_{ir}\right)\\
    &=: \mathcal{K}_{11} + \mathcal{K}_{12}.
\end{aligned}
\]
For $\mathcal{K}_{12}$, we can proceed similarly as before to obtain
\[\mathcal{K}_{12}\le  2\psi(\delta) \Lambda(0)|v_{ir}|.\]

\noindent For $\mathcal{K}_{11}$, since $a_{ik} = 0$ for any $k \neq r, i$ and $a_{rk}=0$ for any $k \neq r$, we have
\[ 
\begin{aligned}
    \mathcal{K}_{11} &= \frac{1}{N}\sum_{k \in C_\ell}\left(a_{ik}\psi_{ik}v_{ki}\cdot v_{ir}-a_{rk}\psi_{rk}v_{kr}\cdot v_{ir}\right)\\
    &= -\frac{1}{N}\psi_{ir}|v_{ir}|^2.
\end{aligned}
\]
Thus, we obtain
\[\psi_{ir}|v_{ir}| \le 2N\psi(\delta)\Lambda(0) - N\frac{d}{dt}|v_{ir}|, \quad \mbox{for a.e. } \ t \in (0,t_0). \]
Hence,
\[\int_{0}^t\psi_{ir}|v_{ir}|d\tau\le (2\psi(\delta)t_0+1)N\Lambda(0)t_0:= M, \quad \mbox{for every } \ t \in [0,t_0). \]

\vspace{0.1cm}

\noindent $\bullet$ (Induction step) Assume that for every $j \in C_\ell \setminus \{r\}$ with $d(r,j) = n-1$, the followings hold:
\[\int_0^t\psi_{j,\pi(j)}|v_{j, \pi(j)}| d\tau\le (n-1)M.\]
Choose any $i \in C_\ell\setminus\{r\}$ with $d(r,i) = n$. Note that $d(r,i) = n$ if and only if
\[
(\underbrace{\pi\circ \cdots \circ \pi}_{n \text{ times}})(i) =r.
\]
Then we have
\[\begin{aligned}
    \frac12\frac{d}{dt} |v_{i,\pi(i)}|^2 &=\frac{1}{N}\left(\sum_{k \in C_\ell} + \sum_{k\not \in C_\ell }\right)\left(a_{ik}\psi_{ik}v_{ki}\cdot v_{i,\pi(i)}-a_{\pi(i),k}\psi_{\pi(i),k}v_{k,\pi(i)}\cdot v_{i,\pi(i)}\right)\\
    &=: \mathcal{K}_{21}+\mathcal{K}_{22}
\end{aligned}\]
For $\mathcal{K}_{22}$, we can similarly get
\[\mathcal{K}_{22}\le 2\psi(\delta) \Lambda(0) |v_{i,\pi(i)}|.\]
For $\mathcal{K}_{21}$, we obtain

\[
\begin{aligned}
    \mathcal{K}_{21} &= \frac{1}{N}\sum_{k \in C_\ell}\left(a_{ik}\psi_{ik}v_{ki}\cdot v_{i,\pi(i)}-a_{\pi(i),k}\psi_{\pi(i),k}v_{k,\pi(i)}\cdot v_{i,\pi(i)}\right)\\
    &=-\frac{1}{N}\psi_{i,\pi(i)}|v_{i,\pi(i)}|^2 - \frac{1}{N}\psi_{\pi(i),\pi(\pi(i))}v_{\pi(i),\pi(\pi(i))}\cdot v_{i,\pi(i)}\\
    &\le -\frac{1}{N}\psi_{i,\pi(i)}|v_{i,\pi(i)}|^2 + \frac{1}{N}\psi_{\pi(i),\pi(\pi(i))}|v_{\pi(i),\pi(\pi(i))}|| v_{i,\pi(i)}|
\end{aligned}
\]
Hence, 
\[\psi_{i,\pi(i)}|v_{i,\pi(i)}| \le 2N\psi(\delta)\Lambda(0)+\psi_{\pi(i),\pi(\pi(i))}|v_{\pi(i),\pi(\pi(i))}| - N\frac{d}{dt}|v_{i,\pi(i)}|.\]
for a.e. $t \in (0,t_0)$ and by integrating the above relation and using the induction hypothesis,
\[\int_0^t \psi_{i,\pi(i)}|v_{i,\pi(i)}| d\tau \le nM\]
for all $t\in [0,t_0)$ and hence the induction is complete. 

Since $d(r,i) \le N-1$ for any $i \in C_\ell\setminus \{r\}$, we can therefore attain

\[\begin{aligned}
|\mathcal{P}_{C_\ell}(t)| &\le |\mathcal{P}_{C_\ell}(0)| + \sum_{i\in C_\ell\setminus\{r\}}\int_0^t\psi_{i,\pi(i)}(\tau)|v_{i,\pi(i)}(\tau)|d\tau\\
&\le |\mathcal{P}_{C_\ell}(0)| + \sum_{n=1}^{N-1}\sum_{\substack{i\in C_\ell\setminus\{r\} \\ d(r,i) = n }}\int_0^t\psi_{i,\pi(i)}(\tau)|v_{i,\pi(i)}(\tau)|d\tau\\
&\le |\mathcal{P}_{C_\ell}(0)| +  \sum_{n=1}^{N-1}\sum_{\substack{i\in C_\ell\setminus\{r\} \\ d(r,i) = n }} nM\\
&\le |\mathcal{P}_{C_\ell}(0)| +  N^2M,
\end{aligned}\]
for all $t \in [0,t_0)$, which contradicts $\mathcal{P}_{C_\ell}(t) \to -\infty$ as $t \nearrow t_0$. This completes the proof of Theorem \ref{thm:glob_ex} under condition (2). 

\vspace{0.2cm}

\subsection{Proof of Theorem \ref{thm:glob_ex} under condition (3)}
Now, we assume that the network topology is undirected. First, we introduce the concept of the Laplacian matrix.

%

\begin{definition}
    Let $M = (m_{ij})_{ij}$ be an $N\times N$ matrix. The \textbf{Laplacian matrix} $L$ of $M$ is defined as
    \[L= D-M,\]
    where $D$ is a diagonal $N\times N$ matrix such that $D_{ii}\displaystyle= \sum_{j=1}^n m_{ij}$ for $i = 1,...,  n$.
\end{definition}

\begin{remark}
 For $\mathds{1}=(1,\cdots, 1)^T \in \R^N$,   every Laplacian matrix $L$ satisfies $L\mathds{1} = 0$. Let $G = (V,E)$ be a digraph with $|V|=N$. Then we can consider the Laplacian matrix $L$ of the $0-1$ adjacency matrix $A = (a_{ij})_{ij}$ of $G$, which is called the (Combinatorial) Laplacian of the graph $G$. If $G$ is an undirected graph, the Laplacian $L$ of $G$ is symmetric and positive semidefinite. Indeed, 
\[
        (Lx)\cdot x =\displaystyle\sum_{i,j=1}^N a_{ij}x_i(x_i-x_j) = \frac{1}{2}\sum_{i,j=1}^N a_{ij}(x_i-x_j)^2\ge 0 \quad\mbox{for all} \quad x \in \mathbb{R}^N. 
\]
    Thus, the eigenvalues $\lambda_i$ of $L$ can be ordered as $0=\lambda_1\le \lambda_2\le \lambda_3\le \cdots \le \lambda_N$ and $L$ admits an orthonormal eigenbasis $\{u_i\}_{i=1}^N$ with $u_1 = \mathds{1}/\sqrt{N}$.
\end{remark}

Then we have the following:

\begin{proposition}\label{prop:lpm_eig}
    Let $G = (V,E)$ be a weakly connected, undirected graph with $|V|=N$. Then we have $\lambda_2 >0$ and
\[
	\frac{\lambda_2}{N}\sum_{i,j=1}^N (x_i-x_j)^2\le \sum_{i,j=1}^N a_{ij}(x_i-x_j)^2 \quad \mbox{for all} \quad x\in \mathbb{R}^N.
\]

\end{proposition}
\begin{proof}
See Corollary 13.4.2 and Lemma 13.1.5 of \cite{Godsil2001}.
\end{proof}

\vspace{0.2cm}

\begin{proof}[Proof of Theorem \ref{thm:glob_ex} under condition (3)] 
As before, we assume for a contradiction that a collision occurs in a finite time and similarly set $t_0$, $C_\ell$, and $\delta>0$.

We define $\|X\|_{C_\ell} $ and $ \|V\|_{C_\ell}$ as follows: 
\begin{center}
    $\displaystyle \|X\|_{C_\ell}(t) := \sqrt{\sum_{i,j \in C_\ell} |x_{ij}(t)|^2}$ \;\;and\;\; $\displaystyle \|V\|_{C_\ell}(t) := \sqrt{\sum_{i,j \in C_\ell} |v_{ij}(t)|^2}.$ 
\end{center}
Then $\displaystyle \|X\|_{C_\ell}(t) \to 0$ as $t \nearrow t_0$. We will prove that the claim in Theorem 1 holds for new $ \|X\|_{C_\ell}$ and $\|V\|_{C_\ell}$. For $t\in(0,t_0)$ at which $ \|X\|_{C_\ell}$ and $\|V\|_{C_\ell}$ are differentiable,

\begin{eqnarray*}
    \left|\frac{d}{dt}  \|X\|_{C_\ell}\right| \le  \|V\|_{C_\ell},
\end{eqnarray*}
and
\[
\begin{aligned}
\frac12\frac{d}{dt}  \|V\|_{C_\ell}^2 &= \sum_{i,j\in C_\ell}v_{ij} \cdot \left[\frac{1}{N}\sum_{k=1}^N a_{ik}\psi_{ik}v_{ki} - \frac{1}{N}\sum_{k=1}^N a_{jk}\psi_{jk}v_{kj}\right]\\
     &= \frac{1}{N}\sum_{i,j\in C_\ell} \sum_{k = 1}^N\left(a_{ik}\psi_{ik}v_{ki}\cdot v_{ij}-a_{jk}\psi_{jk}v_{kj}\cdot v_{ij}\right)\\
    &=\frac{1}{N}\sum_{i,j\in C_\ell}\left(\sum_{k \in C_\ell} + \sum_{k\not \in C_\ell }\right)\left(a_{ik}\psi_{ik}v_{ki}\cdot v_{ij}-a_{jk}\psi_{jk}v_{kj}\cdot v_{ij}\right)\\
    &=: \mathcal{L}_1+\mathcal{L}_2.
\end{aligned}
\]
We separately estimate $\mathcal{L}_1$ as follows:

\[
\begin{aligned}
\mathcal{L}_1 &= \frac{1}{N}\sum_{i,j,k\in C_\ell}a_{ik}\psi_{ik}v_{ki}\cdot v_{ij}-\frac{1}{N}\sum_{i,j,k\in C_\ell}a_{jk}\psi_{jk}v_{kj}\cdot v_{ij} \\
&=: \mathcal{L}_{11}+ \mathcal{L}_{12}.
\end{aligned}\]
For $\mathcal{L}_{11}$, we substitute indices $k$ and $i$ to obtain
\[
\begin{aligned}
  \mathcal{L}_{11} &=\frac{1}{2N}\sum_{i,j,k\in C_\ell}a_{ik}\psi_{ik}v_{ki}\cdot v_{ij} +\frac{1}{2N}\sum_{i,j,k\in C_\ell}a_{ik}\psi_{ik}v_{ik}\cdot v_{kj}\\
    &=\frac{1}{2N}\sum_{i,j,k\in C_\ell}a_{ik}\psi_{ik}v_{ki}\cdot v_{ik}\\
    &=-\frac{|C_\ell|}{2N}\sum_{i,j\in C_\ell}a_{ij}\psi_{ij}|v_{ij}|^2,
\end{aligned}
\]
and similarly,

\[\mathcal{L}_{12}= -\frac{|C_\ell|}{2N}\sum_{i,j\in C_\ell}a_{ij}\psi_{ij}|v_{ij}|^2.\]
Since we have $|x_{ij}|\le  \|X\|_{C_\ell}$ for every $i \neq j$ in $C_\ell$, 
\[\mathcal{L}_1\le -\frac{|C_\ell|}{N}\psi( \|X\|_{C_\ell})\sum_{i,j \in C_\ell}a_{ij}|v_{ij}|^2.\]
For each $i \in C_\ell$, we write $v_i = (v_i^1,v_i^2,..., v_i^d) \in \mathbb{R}^d$. Then since $G[C_\ell]$ is weakly connected, we use Proposition \ref{prop:lpm_eig} to get

\begin{eqnarray*}
    \sum_{i,j\in C_\ell}a_{ij}|v_{ij}|^2 &=& \sum_{i,j\in C_\ell}a_{ij}\sum_{k=1}^d(v_i^k-v_j^k)^2\\
    &=&\sum_{k=1}^d\sum_{i,j\in C_\ell}a_{ij}(v_i^k-v_j^k)^2\\
    &\ge&\frac{\lambda^*}{|C_\ell|} \sum_{k=1}^d\sum_{i,j\in C_\ell}(v_i^k-v_j^k)^2\\
    &=&\frac{\lambda^*}{|C_\ell|}\sum_{i,j\in C_\ell} |v_{ij}|^2
\end{eqnarray*}
where $\lambda^*>0$ is the first positive eigenvalue of the Laplacian of $G[C_\ell]$. Thus, 

\[\mathcal{L}_1\le -\frac{\lambda^*}{N}\psi( \|X\|_{C_\ell}) \|V\|_{C_\ell}^2.\]
For $\mathcal{L}_2$, by using the same argument as in the proof of Theorem 3.1 under condition (1),

\[\mathcal{L}_2\le 2\psi(\delta)\Lambda(0)\sum_{i,j \in C_\ell}|v_{ij}|\le 2\psi(\delta)\Lambda(0)\sqrt{\sum_{i,j \in C_\ell}1}\sqrt{\sum_{i,j \in C_\ell}|v_{ij}|^2}\le 2\psi(\delta)N\Lambda(0) \|V\|_{C_\ell}.\]

\noindent Collecting all the above estimates, we have
\[
\frac{d}{dt} \|V\|_{C_\ell} \le -\frac{\lambda^*}{N}\psi( \|X\|_{C_\ell}) \|V\|_{C_\ell} + 2\psi(\delta)N\Lambda(0).\]
Then we again use Lemma \ref{lem:ineq_ex} and derive a contradiction to obtain the desires result.

\end{proof}

\subsection{One-dimensional case; Proof of Theorem \ref{thm:1d_glob}}

%

Next, we move on to the one-dimensional case and prove the first assertion of Theorem \ref{thm:1d_glob}.

\subsubsection{Proof of Theorem \ref{thm:1d_glob} (i) }
We let $\alpha \ge 1$ and assume that every connected induced subgraph of $G$ has at most one source.

 As in the proof of Theorem \ref{thm:glob_ex}, we assume for a contradiction that a collision occurs in a finite time. We similarly set $t_0$, $C_\ell$, and $\delta>0$ as before. First, we choose indices $m\neq M$ in $C_\ell$ satisfying
  \[ 
 x_m(0) = \min_{i\in C_\ell }x_i(0), \quad \mbox{and} \quad x_M(0)= \max_{i \in C_\ell}x_i(0).
 \]
Since $t=t_0$ is the first collision time, we find that
\[
x_{k}(t) > x_m(t) \quad\mbox{for}\quad k\in C_\ell\setminus \{m\}  \quad \mbox{with} \quad a_{mk}=1
\]
and 
\[
x_{k}(t) < x_M(t) \quad \mbox{for} \quad k\in C_\ell \setminus \{M\} \quad \mbox{with} \quad a_{Mk}=1
\]\\
for all $t\in[0,t_0)$.

Define the sets of indices 
\[
\begin{aligned}
I_m &:= \{k\in C_\ell\setminus\{m\} \mid a_{mk}=1\}, \quad \mbox{and}\\
 I_M &:= \{k\in C_\ell\setminus\{M\} \mid a_{Mk}=1\}.
 \end{aligned}
 \] 
Since $G[C_\ell]$ has at most one source. we have $I_m\cup I_M \neq \emptyset$.  Indeed, if $I_m \cup I_M = \emptyset$, then $I_m = I_M = \emptyset$, so $m$ and $M$ are distinct sources of $G[C_\ell]$, leading to a contradiction.

If $I_m \neq \emptyset$, we get
\[\begin{aligned}
    \dot v_{m} &= \frac{1}{N}\sum_{k=1}^Na_{mk}\frac{v_{km}}{|x_{mk}|^\alpha}\\
    &=\frac{1}{N}\sum_{k\in I_m}(x_{km})^{-\alpha}v_{km} + \frac{1}{N} \sum_{k \notin C_\ell} a_{mk}\frac{v_{km}}{|x_{mk}|^\alpha}\\
    &= \frac{1}{N}\sum_{k\in I_m}\frac{d}{dt}\Psi(x_{km}) + \frac{1}{N}\sum_{k \not \in C_\ell}a_{mk}\frac{v_{km}}{|x_{mk}|^\alpha},
\end{aligned}\]
and hence,
\[
\begin{aligned}
\sum_{k\in I_m}\Psi(x_{km}(t))&=&N(v_m(t)-v_m(0)) + \sum_{k\in I_m}\Psi(x_{km}(0)) - \sum_{k\not \in C_\ell}\int_{0}^ta_{mk}\frac{v_{km}}{|x_{mk}|^\alpha}d\tau
\end{aligned}
\]
for all $t \in [0,t_0)$.
Here, note that 
\bq\label{eq:max_sp}
P(t) := \max_{i \in V}|v_i(t)| \le P(0).
\eq  
Indeed, for a.e. $t \in (0, t_0)$, we can find $1 \le i \le N$ such that $P(t) = |v_i(t)|$ and hence,
\[
\frac12\frac{d}{dt}P^2 =\frac12\frac{d}{dt}|v_i|^2 = \frac{1}{N}\sum_{j=1}^Na_{ij}\psi_{ij}v_{i}\cdot(v_j-v_i)\le \frac{1}{N}\sum_{j=1}^Na_{ij}\psi_{ij}|v_i|(|v_j|-|v_i|)\le 0,
\]
which implies \eqref{eq:max_sp}. Thus, for all $t \in [0,t_0)$,
\bq\label{eq:sum_min}
 \left|\sum_{k\in I_m}\Psi(x_{km}(t))\right| \le 2NP(0) + \sum_{k\in I_m}|\Psi(x_{km}(0)) |+\frac{N \Lambda(0)t_0}{\delta^\alpha}.
\eq

Similarly, if $I_M\neq \emptyset$,
\[\begin{aligned}
    \dot v_{M} &= \frac{1}{N}\sum_{k=1}^Na_{Mk}\frac{v_{kM}}{|x_{Mk}|^\alpha}\\
    &=-\frac{1}{N}\sum_{k\in I_M}(x_{Mk})^{-\alpha}v_{Mk} + \frac{1}{N} \sum_{k \not\in C_\ell} a_{Mk}\frac{v_{kM}}{|x_{Mk}|^\alpha}\\
    &= -\frac{1}{N}\sum_{k\in I_M}\frac{d}{dt}\Psi(x_{Mk}) + \frac{1}{N}\sum_{k \not \in C_\ell}a_{Mk}\frac{v_{kM}}{|x_{Mk}|^\alpha},
\end{aligned}\]
and hence,

\bq\label{eq:sum_max}
    \left|\sum_{k\in I_M}\Psi(x_{Mk}(t))\right| \le 2NP(0) + \sum_{k\in I_M}|\Psi(x_{Mk}(0)) |+\frac{N \Lambda(0)t_0}{\delta^\alpha}.
\eq
for all $t \in [0,t_0)$. Since $\Psi(r) \to -\infty$ as $r \to 0^+$ and $I_m \cup I_M \neq \emptyset$, at least one of the left-hand side of \eqref{eq:sum_min} and \eqref{eq:sum_max} must diverge to infinity, which leads to a contradiction.

\subsubsection{Proof of Theorem \ref{thm:1d_glob} (ii)} Now, we let $\alpha>0$ By assumption, we have a subset of vertices $A\subset V$ such that $G[A]$ is weakly connected and has at least two sources. Then $G[A]$ has at least three vertices, i.e. $|A|\ge 3$. Let $s_1$ and $s_2$ be two distinct sources of $G[A]$. Since $G[A]$ is weakly connected and $s_1,s_2$ are sources, the length of the shortest path from $s_1$ to $s_2$ in $G^*[A]$ must be at least $2$.

Then we can find vertices $i_1, ..., i_k \in A$ satisfying $s_1 \to i_1 \to \cdots  \to i_k \to s_2$ forms the shortest path from $s_1$ to $s_2$ in $G^*[A]$. If necessary, we relabel vertices so that vertices $1$ and $k+2$ are sources of $G[A]$ and $1 \to 2 \to \cdots \to k+2$ is the shortest path from $1$ to $k+2$. Now, construct the initial data as follows:
\begin{itemize}
    \item $x_1(0)=-1$, $v_1(0)=10$,
    \item $x_{k+2}(0)=1$, $v_{k+2}(0)=-10$,
    \item For $i \in \{2,...,k+1\}$, $x_i(0) = -1 \displaystyle+\frac{2(i-1)}{k+1}$, $v_i(0)=0$,
    \item For $j \in \{k+3,...N\}$, $x_j(0)=10j$, $v_j(0)=0$. 
\end{itemize}
Now, assume for a contradiction that a global-in-time smooth solution to \eqref{main_sys} with the initial data above exists. Let $P(t) := \displaystyle\max_{i \in V}|v_i(t)|$. From \eqref{eq:max_sp}, we have $P(t) \le P(0)$ and hence, $|v_i(t)|\le 10$ for every $i=1, \dots, N$. Then it follows that, for each $p=1,...,k+2$ and $q = k+3,...,N$, 
    \[|x_p-x_q|\ge 25\text{ on $[0,1/5]$}.\]
Since $\{1,...,k+2\}\subseteq A$ and vertices $1$ and $k+2$ are sources of $G[A]$, we have $a_{1j} = a_{k+2, j}=0$ for every $j \in A$. Thus, for every $t \in [0,1/5]$,

\[\begin{aligned}
|\dot v_1| &= \left|\frac{1}{N}\sum_{j=1}^N a_{1j}\frac{v_{j1}}{|x_{1j}|^\alpha}\right| \\
&= \left|\frac{1}{N}\sum_{j\not \in A} a_{1j}\frac{v_{j1}}{|x_{1j}|^\alpha}\right|\\
&\le\frac{1}{N}\sum_{j\not \in A} a_{1j}\frac{|v_{j1}|}{|x_{1j}|^\alpha}\le  \frac{1}{N}\sum_{j\not \in A}a_{1j}\frac{20}{25^\alpha} \le \frac{20}{25^\alpha}.
\end{aligned}\]
Similarly, we can also get $|\dot v_{k+2}|\le \displaystyle\frac{20}{25^\alpha}$ and this leads to
\[
|(\dot v_{k+2}-\dot v_1)(t)|\displaystyle\le \frac{40}{25^\alpha}, \quad \forall t \in [0, 1/5].
\] 
Thus, for every $t \in [0,1/5]$,

\[\begin{aligned}
    x_{k+2}(t)-x_1(t) &= x_{k+2}(0)-x_1(0)+\int_{0}^t [v_{k+2}(\tau)- v_1(\tau)] d\tau\\
     &= x_{k+2}(0)-x_1(0)+t(v_{k+2}(0)- v_1(0) )+\int_{0}^t\int_{0}^\tau [\dot v_{k+2}(\sigma)-\dot v_1(\sigma)]d\sigma d\tau\\
     &= 2-20t +\int_{0}^t\int_{0}^\tau [\dot v_{k+2}(\sigma)-\dot v_1(\sigma)]d\sigma d\tau\\
     &\le 2-20t +\frac{1}{2}\cdot \frac{40}{25^\alpha} t^2.
\end{aligned}\]

This implies $x_{k+2}(1/5)-x_{1}(1/5) \le -2 +\displaystyle\frac{20}{25^{\alpha+1}}<0$. If we write $u_i(t) := (x_{i+1} - x_i)(t)$ for each $i=1,\dots, k+1$, then 
\[
\sum_{i=1}^{k+1} u_i(0) = (x_{k+2} - x_1)(0) >0, \quad \mbox{and} \quad \sum_{i=1}^{k+1} u_i(1/5) = (x_{k+2} - x_1)(1/5) <0.
\]
Since each $u_i$ is continuous, the intermediate value theorem implies that there exists $t_0\in (0,1/5)$ satisfying
\[
\sum_{i=1}^{k+1}u_i(t_0) =0,
\]
Here, either $u_i(t_0) = 0$ for all $i=1,\dots, k+1$ or $u_i(t_0)<0$ for some $i$. Since $u_i(0)>0$ for all $i=1,\dots, k+1$, by the intermediate value theorem if necessary, we can find an index $j$ and  $t^* \in (0,1/5)$ satisfying $u_j(t^*) = 0$. As there is a directed path $j \to (j+1)$ or $(j+1) \to j$ by assumption, this implies agents $j$ and $j+1$ collide at $t=t^*$ according to Definition \ref{def:col}. Therefore, it contradicts the global-in-time existence of solutions and this completes the proof.

\vspace{0.2cm}

%
%

%
%
%
%

\section{Emergence of mono-cluster flocking}\label{sec:flock}
In this section, we present the proof of Theorems \ref{thm:LT1} and \ref{thm:LT2}. First,

\subsection{Proof of Theorem \ref{thm:LT1}}
First, we review the property of stochastic matrices, which are nonnegative matrices whose row sums are all equal to 1. 
\begin{proposition}\cite{DHK19}\label{prop:eg_coeff}
    Assume that an $N\times N$ matrix $M = (m_{ij})$ is stochastic, and let $W = (w_i^j)$ and $Z = (z_{i}^j)$ be $N\times d$ matrices such that
    \[W = MZ\]
    Then, we have
    \[\max_{i,j} | w_i-w_j| \le (1-\mu(M))\max_{i,j} |z_i-z_j|\]
    where $w_i := (w_i^1, \cdots,w_{i}^d),z_i:= (z_i^1,\cdots, z_i^d), i=1,\cdots, N$ and $\mu$ denotes the ergodicity coefficient (see Remark \ref{rm:scrb}).
\end{proposition}
Then we are ready to provide the proof.

\begin{proof}[Proof of Theorem \ref{thm:LT1}]
Since the initial data is non-collisional, we can find $T_0>0$ such that system \eqref{main_sys} admits a unique, smooth solution on the time interval $[0,T_0)$. Without loss of generality, we may let $T_0$ the maximal life-span of system \eqref{main_sys}, i.e.
\[
T_0 := \sup \{ T>0 \ : \ \mbox{system \eqref{main_sys} admits a unique smooth solution on } \ [0,T)\}.
\]
 We prove that the solution exists globally in time, i.e. $T_0 = \infty$ and exhibits the mono-cluster flocking. For this, it sufiices to show the following claim:
\begin{claim}
    For every $m \in \mathbb{N}\cup\{0\}$ with $\hat{h}m \le T_0 $,
\bq\label{eq:vd_decay}
\Lambda(\hat{h}m) \le (1-\delta)^m \Lambda(0),
\eq
 where $\displaystyle \hat{h} := \frac{h}{\psi_\eta}$ and $h$ is the shortest height of all spanning trees of $G$.
\end{claim}

Once we show the claim holds, we can deduce that $T_0 =\infty$. Indeed, if we assume for a contradiction that $T_0<\infty$, we set $m^* := \lfloor T_0 / \hat{h}\rfloor$. Then by \eqref{eq:vd_decay} and Proposition \ref{prop:dm_bd}, one obtains 
\[\Lambda(t) \le \Lambda(\hat{h}k) \le (1-\delta)^k\Lambda(0).\]
for $0\le k\le m^*$ and $\hat{h}k \le t < \min\{\hat{h}(k+1),T_0\}$. Now, choose any pair $i\neq j$ in $V_G$ with $(j,i) \in E_G$ or $(i, j) \in E_G$. Then for every $0\le k\le m^*$ and $\hat{h}k \le t < \min\{\hat{h}(k+1),T_0\}$, we can get
\[
\begin{aligned}
    |x_{ij}(t)|&\ge |x_{ij}(0)| - \int_{0}^t \Lambda(s)ds\\
    &\ge |x_{ij}(0)| - \sum_{k=0}^{m^*-1}\int_{\hat{h}k}^{\hat{h}(k+1)} \Lambda(s)ds - \int_{\hat{h}m^*}^{T_0} \Lambda(s)ds\\
    &\ge |x_{ij}(0)| -\hat{h}\Lambda(0) \sum_{k=0}^{m^*}  (1-\delta)^k\\
    &\ge |x_{ij}(0)| -\hat{h}\Lambda(0) \sum_{k=0}^\infty  (1-\delta)^k\\
    &= |x_{ij}(0)| -\frac{\hat{h}}{\delta}\Lambda(0) \ge \gamma(0) - \eta >0.
\end{aligned}
\]
Thus, we have
\[
\inf_{t\in[0,T_0)} \gamma(t)\ge \gamma(0)-\eta >0,
\]
and this suggests that the unique smooth solution can be extended beyond $T_0$, contradicting the maximality of $T_0$.\\

\vspace{0.2cm}
Then we show the claim holds. It is vacuously true when $m =0$.  Assume that it holds for every $0\le k\le m$ with $\hat{h}m\le T_0$. Since it becomes again vacuously true when $\hat{h}(m+1)>T_0$, we may assume $\hat{h}(m+1)\le T_0$. Then the same argument implies that for $i\neq j$ in $V_G$ with $(j,i) \in E_G$ or $(i,j)\in E_G$, we obtain 
\bq\label{eq:bd_psi}
\psi_\rho \le \psi_{ij}(t) \le \psi_\eta
\eq
for $t\in [0,\hat{h}(m+1)]$. 

 Next, for system \eqref{main_sys} with $X = (x_1,x_2,...,x_N)$ and $V = (v_1,v_2,...,v_N)$, we rewrite it as 
\bq\label{eq:mat_sys}
\begin{aligned}
\frac{d}{dt}X &= V,\\
\frac{d}{dt}V &= -\overline{L}V,
\end{aligned}
\eq
where  $\overline{L} := D-R$ is an $N\times N$ matrix and $D$ and $R$ are defined as
\[
R_{ij} := \lt\{\begin{array}{lcl}\displaystyle\frac1N a_{ij} \psi_{ij} & \mbox{if} & i \neq j,\\
& & \\
\displaystyle\frac1N \psi_\eta & \mbox{if} & i=j, \end{array} \rt.  \quad D_{ij} := \lt\{\begin{array}{lcl}0 & \mbox{if} & i \neq j,\\
& & \\
\displaystyle\sum_{k=1}^N R_{ik} & \mbox{if} & i=j, \end{array} \rt.
\]
so that $\overline{L}$ is the Laplacian matrix of $R$. 

Furthermore, we recast the second equation in \eqref{eq:mat_sys} as
\[
 \frac{d}{dt}V = -\overline{L} V \quad \Longleftrightarrow \quad \frac{d}{dt}(e^{\psi_{\eta}t} V)= (\psi_{\eta}\mathbb{I}_N - \overline{L})(e^{\psi_{\eta}t} V).
\]
Note that the matrix $(\psi_{\eta}\mathbb{I}_N - \overline{L})$ is nonnegative on the time interval $[0,\hat{h}(m+1)]$ due to \eqref{eq:bd_psi}. Here we deduce from \eqref{eq:stat_tran} that
\[V(t) = e^{-\psi_\eta (t-t_0)} \Phi(t,t_0)V(t_0),\]
for any $0 \le t_0 \le t \le \hat{h}(m+1)$, where
\[
\begin{aligned}
    \Phi(t,t_0) &= \mathbb{I}_N + \int_{t_0}^t
    (\psi_\eta \mathbb{I}_N - \overline{L}(s))ds \\
&\quad+ \sum_{k=1}^\infty \int_{t_0}^t \int_{t_0}^{s_0} \cdots \int_{t_0}^{s_{k}} (\psi_\eta \mathbb{I}_N - \overline{L}(s_0))(\psi_\eta \mathbb{I}_N - \overline{L}(s_1))\cdots (\psi_\eta\mathbb{I}_N - \overline{L}(s_{k+1}))ds_{k+1}ds_{k}\cdots ds_0
\end{aligned}
\]
Here, the row sums of the matrix $\tilde{\Phi}(t,t_0) := e^{-\psi_\eta (t-t_0)} \Phi(t,t_0)$ are equal to 1. This allows us to use Proposition \ref{prop:eg_coeff} so that
\[
\Lambda(\hat{h}(m+1))\le \lt[1-\mu(\tilde{\Phi}(\hat{h}(m+1),\hat{h}m))\rt]\Lambda(\hat{h}m),
\]
and one uses $(\psi_\eta \mathbb{I}_N - \overline{L}) \ge R\ge \displaystyle\frac{\psi_\rho}{N}A$, where $A$ is the 0-1 adjacency matrix of $G$, to obtain
\[
\begin{aligned}
    \Phi(t,t_0) &\ge \int_{t_0}^t \int_{t_0}^{s_0}\cdots \int_{t_0}^{s_{h-2}} (\psi_\eta \mathbb{I}_N - \overline{L}(s_0))(\psi_\eta \mathbb{I}_N - \overline{L}(s_1))\cdots (\psi_\eta \mathbb{I}_N - \overline{L}(s_{h-1}))ds_{h-1}ds_{h}\cdots ds_0\\
    &\ge  \left(\frac{\psi_\rho}{N}\right)^hA^h\int_{t_0}^t \int_{t_0}^{s_0}\cdots \int_{t_0}^{s_{h-2}} ds_{h-1}ds_{h}\cdots ds_0\\
    &=\frac{(t-t_0)^h}{h!}\left(\frac{\psi_\rho}{N}\right)^hA^h.
\end{aligned}
\]
where $h$ denotes the shortest height of all spanning trees of $G$. Then by applying Remark \ref{rm:scrb}, we have
\[\mu(\tilde{\Phi}(\hat{h}(m+1),\hat{h}m))\ge\mu\left(\frac{1}{h!}\left(\frac{h\psi_\rho}{Ne \psi_\eta}\right)^hA^h\right)\ge \frac{1}{h!}\left(\frac{h\psi_\rho}{Ne \psi_\eta}\right)^h = \delta.\]

Thus,
 \[\Lambda(\hat{h}(m+1))\le (1-\delta)\Lambda(\hat{h}m)\le (1-\delta)^{m+1}\Lambda(0),\]
and this completes the induction. Therefore, for $m \in \mathbb{N}\cup\{0\}$, 
\[
\Lambda(\hat{h}m) \le (1-\delta)^m \Lambda(0)
\]
and this implies 
\[
\Gamma(t)\le \Gamma(0) + \rho \quad \mbox{and} \quad \gamma(t)\ge \gamma(0) - \eta,
\] 
for all $t\ge 0$. 
 \end{proof}

\vspace{0.2cm}

\subsection{Proof of Theorem \ref{thm:LT2}}
We proceed to the case where the network topology is scrambling. First,  we let $T_0$ be the maximal life-span of system \eqref{main_sys} for the given initial data $(X(0),V(0))$. Note that Theorem \ref{thm:glob_ex} implies  $T_0 = \infty$ under conditions (ii) and (iii). Now, for $t\in(0,T_0)$ at which $ \Lambda(t)$ is differentiable, we may let $ \Lambda(t)=|v_{ij}|$ for some $i,j \in V_G$. Then

\[
\begin{aligned}
    \frac12\frac{d}{dt} \Lambda^2(t) &= \frac12\frac{d}{dt}|v_{ij}|^2 = \frac{1}{N}\sum_{k=1}^N(a_{ik} \psi_{ik}v_{ki} \cdot v_{ij} - a_{jk}\psi_{jk}v_{kj}\cdot v_{ij})\\
    &\le -\frac{1}{N}\psi(\Gamma(t))\sum_{k=1}^N(a_{ik}v_{ik} \cdot v_{ij} + a_{jk}v_{kj}\cdot v_{ij})\\
    &\le -\frac{1}{N}\psi(\Gamma(t))\sum_{k \to i,j}(v_{ik} \cdot v_{ij} + v_{kj}\cdot v_{ij})\\
    &= -\frac{1}{N}\psi(\Gamma(t))\Lambda^2(t)\sum_{k \to i,j}1\\
    &\le -\frac{S(G)}{N}\psi(\Gamma(t))\Lambda^2(t).
\end{aligned}
\]
Thus, we obtain
\bq\label{eq:est_dv_scr}
\frac{d}{dt}\Lambda(t) \le -\frac{S(G)}{N}\psi(\Gamma(t))\Lambda(t).
\eq
Motivated by \cite{ACHL12, HL09}, we define a Lyapunov functional as
 
\[
\mathcal{L}(t) := \frac{N}{S(G)}\Lambda(t) + \int_{\Gamma(0)}^{\Gamma(t)} \psi(s)\,ds.
\]
Then we have
\[\displaystyle\frac{d}{dt}\mathcal{L}(t)=\frac{N}{S(G)}\frac{d}{dt}\Lambda(t) + \psi(\Gamma(t)) \frac{d}{dt}\Gamma(t) \le -\psi(\Gamma(t))\Lambda(t) + \psi(\Gamma(t))\Lambda(t)=0.\]
Hence,
\[
\int_{\Gamma(0)}^{\Gamma(t)} \psi (s)\,ds\le \frac{N}{S(G)} \Lambda(0).
\]
Under condition (i), non-integrability of $\psi$ implies the existence of $\Gamma_*$ satisfying
\[
\int_{\Gamma(0)}^{\Gamma_*}\psi(s)\,ds = \frac{N}{S(G)} \Lambda(0),
\]
which coincides with that in the statement of Theorem. This gives
\[
\Gamma \le \Gamma_* \quad \mbox{on} \quad [0,T_0),
\]
and hence,
\[
\Lambda(t)\le \Lambda(0)\exp\left(-\frac{S(G)}{N\Gamma_*^\alpha}t\right) \quad \forall t \in [0,T_0).\]
Moreover, we can deduce that, for $t\in [0,T_0)$,
\[
\begin{aligned}
    \gamma(t) &\ge \gamma(0)- \int_0^t \Lambda(s)\,ds\\
    &\ge\gamma(0) - \frac{N}{S(G)\psi(\Gamma_*)}\Lambda(0)\left(1-\exp\left(-\frac{S(G)}{N}\psi(\Gamma_*)t\right)\right)\\
    &\ge\gamma(0) - \frac{N}{S(G)\psi(\Gamma_*)}\Lambda(0)>0
\end{aligned}
\]
This suggests that $T_0=\infty$ and also the mono-cluster flocking estimates. The results under conditions (ii) and (iii) readily follow from the above arguments and this completes the proof.


\vspace{0.4cm}

\noindent Finally, we can apply similar arguments to the cases where the network topology is undirected.
\begin{corollary}\label{cor:und_LT}
Let $\alpha>0$ and $G= (V_G, E_G)$ be the undirected, weakly connected graph representing the network topology in \eqref{main_sys}. Suppose that the initial data $(X(0), V(0))$ are non-collisional,  i.e. $\gamma(0)>0$, and one of the following conditions holds:
    \begin{itemize}
        \item[(i)] $\alpha \in (0,1)$, and 
        \[\|\tilde{V}(0)\| := \lt(\sum_{i=1}^N (v_i -v_c(0))^2 \rt)^{1/2}<\frac{\lambda_2\gamma(0)}{\sqrt{2}N \tilde{\Gamma}_*^\alpha}\]
        \item[(ii)] $\alpha>1$, and
        \[\|\tilde{V}(0)\|<\frac{\lambda_2\Gamma(0)^{1-\alpha}}{\sqrt{2}N(\alpha-1)}, \]
        \item[(iii)] $\alpha=1$,
    \end{itemize}
    where $v_c := \frac1N \sum_{i=1}^N v_i$ is the mean velocity, $\lambda_2$ is the first positive eigenvalue of the Laplacian of $G$ and $\tilde{\Gamma}_*$ is defined as
    \[\tilde{\Gamma}_* = \begin{cases}
        \displaystyle\Gamma(0)\exp\left(\frac{\sqrt{2}N\|\tilde{V}(0)\|}{\lambda_2}\right)&,\; \alpha =1\\
        \displaystyle \left(\Gamma(0)^{1-\alpha}+\frac{(1-\alpha)\sqrt{2}N} {\lambda_2}\|\tilde{V}(0)\|\right)^{1/(1-\alpha)}&,\; \alpha \in (0,1) \cup (1, \infty).
    \end{cases}\]
     Then system \eqref{main_sys} admits a global-in-time smooth solution and exhibits mono-cluster flocking:
    \[\Gamma(t) \le \tilde{\Gamma}_*, \quad  \|\tilde{V}(t)\| \le  \|\tilde{V}(0)\|\exp\left(-\frac{\lambda_2}{N\tilde{\Gamma}_*^\alpha}t\right) \quad \mbox{for all} \quad t \ge0.\]

\end{corollary}
\begin{proof}
Let $T_0$ be the maximal life-span of system \eqref{main_sys} for the given initial data $(X(0),V(0))$ as before. When the network topology is given by an undirected graph, the mean velocity is preserved in time:

\[\begin{aligned}
    \frac{d}{dt}\sum_{i=1}^N v_i &= \frac{1}{N} \sum_{i,j=1}^N a_{ij}\psi_{ij} v_{ji} \\
&= \frac{1}{N}\sum_{i,j=1}^N a_{ji}\psi_{ji}v_{ij} \\
&= -\frac{1}{N}\sum_{i,j=1}^N a_{ij}\psi_{ij}v_{ji}= 0.
\end{aligned}\]
So we may assume $v_c(t) \equiv 0$ without loss of generality. Then we use Proposition \ref{prop:lpm_eig} to yield

\[
\begin{aligned}
    \frac12\frac{d}{dt} \|\tilde{V}\|^2 &= \frac{1}{N}\sum_{i,j=1}^N a_{ij} \psi_{ij} v_{ji}\cdot v_i\\
    &= -\frac{1}{2N}\sum_{i,j=1}^N a_{ij} \psi_{ij} |v_{ij}|^2\\
    &\le -\frac{1}{2N}\psi(\Gamma(t))\sum_{i,j=1}^N a_{ij}|v_{ij}|^2\\
    &\le -\frac{\lambda_2}{2N^2}\psi(\Gamma)\sum_{i,j=1}^N|v_{ij}|^2\\
    &= -\frac{\lambda_2 }{N}\psi(\Gamma) \|\tilde{V}\|^2,
\end{aligned}
\]
which implies
\[\frac{d}{dt}\|\tilde{V}\| \le -\frac{\lambda_{2}}{N}\psi(\Gamma)\|\tilde{V}\|.\]
Here, we use $\Lambda \le \sqrt{2}\|\tilde{V}\|$ and set the lyapunov functional as
\[
\mathcal{L} := \frac{\sqrt{2}N}{\lambda_2}\|\tilde{V}\| + \int_{\Gamma(0)}^{\Gamma(t)} \psi(s)\,ds.
\] 

Then we argue similarly to the proof of Theorem \ref{thm:LT2} to obtain the desired results.
\end{proof}


\section*{Acknowledgments}
The work of J. Jung was supported by the National Research Foundation of Korea (NRF) grant funded by the Korean government (MSIT) (RS-2026-25589791).

\appendix
\section{Characterization of the class of digraphs}\label{sec:app.A}
In this appendix, we provide the characterization of the class $\mathcal{S}$ of digraphs on $N$ vertices whose every weakly connected induced subgraph is scrambling and every vertex has a self-loop. We first consider a characterization of $\mathcal{S}$.
\begin{proposition}\label{prop:ch_S1}
    A digraph $G = (V_G,E_G)$ belongs to the class $\mathcal{S}$ if and only if both conditions hold:
    \begin{itemize}
        \item[(1)] $G^*$ is $P_4$-free, i.e. $P_4$ can not be an induced subgraph of $G^*$.
        \item[(2)] For every $\tilde{V} \subset V$ with $G^*[\tilde{V}]$ being a $P_3$, $G[\tilde{V}]$ is scrambling,
    \end{itemize}
    where $G^*$ denotes the underlying graph of $G$ and $P_k$ denotes the path graph $v_1-v_2- \cdots-v_k$ on $k$ vertices.
\end{proposition}
\begin{proof} 
\noindent $(\implies)$ Suppose that $G$ belongs to the class $\mathcal{S}$. If there exists $\tilde{V}\subset V$ such that $G^*[\tilde{V}]$ is a $P_4$, $G^*[\tilde{V}]$ should be of the form 
\[
i_1 - i_2 - i_3 - i_4.
\] 
Then in this graph, there is no $r \in \{i_1, i_2, i_3, i_4\}$ such that both $(r, i_1)$ and $(r, i_4)$ are in the edge set of $G[\tilde{V}]$. So it contradicts the assumption that any weakly connected induced subgraph must be scrambling. Moreover, if $G^*[\tilde{V}]$ is a $P_3$, $G[\tilde{V}]$ is weakly connected, and hence it must be scrambling.\\

\noindent $(\impliedby)$ Now, suppose that $G$ satisfies conditions (1) and (2). Let $H = (V_H,E_H)$ be a connected induced subgraph of $G$. It suffices to show that, for any $x,y \in V_H$ with $(x,y),(y,x) \not \in E_H$, there exists $r \in V_H$ such that $(r,x),(r,y) \in E_H$. The other cases are clear since we assume that $G$ has a self-loop at each of its vertex. By condition (1), the shortest path from $x$ to $y$ in $H^*$ must be a $P_3$, so condition (2) implies the existence of such vertex.
 \end{proof}

To proceed further, we present some concepts related to connected components of digraphs.
\begin{definition}
    Let $G =(V_G,E_G)$ and $H=(V_H,E_H)$ be digraphs. We define $G\cup H$ as a digraph with a vertex set $V = V_G \cup V_H$ and an edge set $E = E_G \cup E_H$. If $V_G$ and $V_H$ are disjoint, we say that $G$ and $H$ are disjoint and $G\cup H$ is a disjoint union of $G$ and $H$.
\end{definition}
Once we recall the notion of a weakly connected component of a digraph from Definition \ref{def:cnt_gp}, we have the following:

\begin{proposition}
    Every digraph $G = (V_G,E_G)$ is the disjoint union of weakly connected components of $G$.
\end{proposition}
\begin{proof} 
Let $G^*$ be the underlying graph of $G$. By identifying each edge with a  $1$-cell, we can regard $G^*$ as the one dimensional CW complex. Since topological connectedness and combinatorial connectedness coincide for such complexes, the topological connected components of $G^*$ are precisely the weakly connected components. Hence, there exists a partition $V_G=V_1\sqcup\cdots\sqcup V_n$ such that each induced subgraph $G^*[V_i]$ is a weakly connected component of $G^*$. Then by definition, there are no directed edges in $G$ between distinct vertex sets $V_i$ and $V_j$. Consequently, $G$ is the disjoint union of the induced subgraphs $G[V_1],\ldots,G[V_n]$, each of which is a weakly connected component of $G$.
\end{proof}

Next, we also need the following concepts:

\begin{definition}
    Let $H = (V_H,E_H)$ be an undirected graph. Then the \textbf{complement graph} of $H$ is an undirected graph $\overline{H} = (V_H,\overline{E}_H)$ defined by $\overline{E}_H = \{(i,j) \in V_H\times V_H\mid (i,j)\not\in E_H \} = (V_H\times V_H) \setminus E_H$. A \textbf{cograph} is defined recursively as follows:
    \begin{enumerate}
        \item  A graph on a single vertex is a cograph.
        \item If $H_1$ and $H_2$ are cographs and disjoint, then $H_1\cup H_2$ is also a cograph.
        \item If $H$ is a cograph, then so is $\overline{H}$. 
    \end{enumerate}
\end{definition}

Now,  we introduce an operation to construct the class $\mathcal{S}$ using semicomplete graphs.
\begin{definition}
    Let $A = (V_A,E_A)$ and $B = (V_B,E_B)$ be disjoint digraphs. 
    \begin{enumerate}
    \item
    We write $\{A \implies B\}$ as the family of digraphs $G=(V_A \cup V_B,E)$ satisfying
    \begin{itemize}
    \item[--]
     For every $a\in V_A$ and $b \in V_B$, $(a,b)\in E$.
     \item[--]
     $G[V_A]$ and $G[V_B]$ are equal to $A$ and $B$, respectively.
     \end{itemize}
     \item
      Moreover, we write $A \iff B$ as the unique element in the intersection of $\{A \implies B\}$ and $\{B \implies A\}$.
      \end{enumerate}
\end{definition}

Using the above operation, we define $\mathcal{T}$ as the smallest class of digraphs on $N$ vertices satisfying the followings:
\begin{enumerate}
\item
Every vertex of a digraph in $\mathcal{T}$ has a self-loop.
    \item
     Every semicomplete graph on $N$ vertices is in $\mathcal{T}$.
    \item
     If $T_1, T_2 \in \mathcal{T}$ are disjoint, their disjoint union also belongs to $\mathcal{T}$,
    \item
     If $W$ is a semicomplete graph, $T \in \mathcal{T}$, and $T$ and $W$ are disjoint, then $\{W \implies T\}\subset \mathcal{T}$.
    \item
     If $T_1, T_2 \in \mathcal{T}$ with $T_1$ and $T_2$ being disjoint, then $(T_1 \iff T_2) \in \mathcal{T}$. 
\end{enumerate}
Note that the operation $\iff$ is associative.
\vspace{0.2cm}

Then we are ready to state the equivalence of $\mathcal{T}$ and $\mathcal{S}$.

\begin{theorem}
    $\mathcal{T} =\mathcal{S}$.
\end{theorem}
\begin{proof} 
The proof consist of two parts; 

\noindent $\bullet$ (Part I: $\mathcal{T}\subseteq \mathcal{S}$) Since $\mathcal{T}$ is defined recursively, we will show that the class $\mathcal{S}$ contains every semicomplete graph and is closed under the aforementioned operations.

\noindent $(1)$ Recall from Remark \ref{rmk:sc_sc} that if $G$ is semicomplete, every induced subgraph of $G$ is scrambling. Thus, any semicomplete graphs belong to $\mathcal{S}$.

\noindent $(2)$ For disjoint $T_1,T_2 \in \mathcal{S}$, every weakly connected induced subgraph of $T_1\cup T_2$ lies entirely in one of $T_1$ and $T_2$. So it is scrambling and hence $T_1\cup T_2\in \mathcal{S}$.

\noindent  $(3)$ Let $W= (V_W, E_W)$ be semicomplete, $T=(V_T, E_T)\in \mathcal{S}$, and $W$ and $T$ are disjoint. Choose any $G \in \{ W \implies T\}$ and consider a weakly connected induced subgraph $H=(V_H, E_H)$ of $G$. If $V_H \subseteq V_W$ or $V_H \subseteq V_ T$, then $H$ is scrambling. Suppose $V_H \cap V_W \neq \emptyset$ and $V_H \cap V_T \neq \emptyset$.  

Choose $x$ and $y$ in $V_H$. If $x$ and $y$ are in $V_W$, then one of them is their common predecessor since every vertex has a self-loop. If $x \in V_W$ and $y \in V_T$, then $x$ is their common predecessor. It is similar when $y \in V_W$ and $x \in V_T$. If $x, y \in V_T$, it suffices to consider the case $(x,y), (y,x) \notin E_H$. If so, any $r \in V_H \cap V_W$ is their common predecessor. Hence $H$ is scrambling and this implies $\{W \implies T\}\subseteq \mathcal{S}$.

\noindent $(4)$ Let $T_1 = (V_{T_1}, E_{T_1}), T_2 = (V_{T_2}, E_{T_2}) \in \mathcal{S}$ be disjoint and put $G := T_1 \iff T_2$. Let $H = (V_H, E_H)$ be any weakly connected induced subgraph of $G$. If $V_H \subseteq V_{T_1}$ or $V_H \subseteq V_{T_2}$, then $H$ is scrambling. Otherwise, take any $x, y \in V_H$. If $x$ and $y$ lie in the different factor, then $x$ is a common predecessor. If they lie in the same factor, any vertex of $H$ in the other factor becomes their common predecessor. Hence $H$ is scrambling and $G \in \mathcal{S}$.

Therefore, we have $\mathcal{T}\subseteq \mathcal{S}.$\\

\noindent $\bullet$ (Part II: $\mathcal{S}\subseteq \mathcal{T}$) Choose any $G = (V_G, E_G) \in \mathcal{S}$. We proceed by induction on $n= |V_G|$.  When $n=1$,  then $G$ is semicomplete and belongs to $\mathcal{T}$.
Assume that every graph in $\mathcal{S}$ with fewer than $n$ vertices belongs to $\mathcal{T}$. Let $G \in \mathcal{S}$ have $n$ vertices. If $G$ is disconnected, we can find a partition $X_1$, $X_2$, ... $X_k$ of $V_G$ such that each $G[X_i]$ is a weakly connected component of $G$. Then since each $G[X_i]$ belongs to $\mathcal{S}$ with $|X_i| <n$ and the class $\mathcal{T}$ is closed under disjoint union, $G \in \mathcal{T}$.

Now, we assume that $G$ is connected. By Proposition \ref{prop:ch_S1}, $G^*$ is $P_4$-free. Then by \cite[Theorem 2]{CLB81}, $G^*$ has the following properties:
\begin{enumerate}
\item
$G^*$ is a cograph.
\item
Any weakly connected subgraph of $G^*$ with at least two vertices has disconnected complement.
\end{enumerate}
We set $C:= \overline{G^*}$ and $V_1, \dots, V_m$ be the partition of $V_G$ such that each $C[V_i]$ is a weakly connected component of $C$. Then $m \ge 2$, each $V_i$ is a nonempty proper subset of $V_G$, and every vertex of $V_i$ is a neighbor of every vertex outside $V_i$ in $G^*$.

Here, we claim that 
\bq\label{eq:cnt_Vi}
|V_i|\ge 2  \implies  (u,v) \in E_G \quad \mbox{for every} \quad u \in V_G\setminus V_i \quad \mbox{and}  \quad v \in V_i.
\eq

For this, we fix $i$, $u\in V\setminus V_i$ and $v \in V_i$. Since $C[V_i]$ is connected and $|V_i|\ge 2$, $v$ has a neighbor $w\in V_i \setminus \{v\}$. This implies $(v,w), (w,v) \notin E_G$. On the other hand, $u$ is a neighbor of $v$ and $w$ in $G^*$ since $u$ belongs to a different component of $C$. Hence, $G^*[\{v, u, w\}]$ is $P_3$ with the middle vertex $u$.

 Then the digraph $G[\{v, u, w\}]$ is a weakly connected induced subgraph of $G$ and hence is scrambling. This implies $(u,v), (u,w) \in E_G$ and proves \eqref{eq:cnt_Vi}.
 
 Now, we define two collections as
 \[
 \mathcal{U} := \bigcup_{|V_i| =1}V_i, \quad \mathcal{W} := \bigcup_{|V_i|\ge 2} V_i.
 \]
The sets in $\mathcal{U} \cup \mathcal{W}$ partition $V_G$. If $\mathcal{U} \neq \emptyset$, then $G[\mathcal{U}]$ is semicomplete. Indeed, any two distinct vertices of $\mathcal{U}$ belong to distinct connected components of $C$ and hence, there must be an edge joining them in $G^*$. 

Suppose $\mathcal{W}\neq \emptyset$ and write the components of $\mathcal{W}$ as $V_{i_1}, \dots, V_{i_\ell}$, i.e.
\[
\mathcal{W} = \bigcup_{k=1}^\ell V_{i_k}.
\]
For each $k=1, \dots, \ell$, since every weakly connected induced subgroup of $G[V_{i_k}]$ is that of $G$, $G[V_{i_k}]$ lies in $\mathcal{S}$.  Hence $G[V_{i_k}] \in \mathcal{T}$ by induction.  Furthermore,  for $V_{i_\alpha}$ and $V_{i_\beta}$ with $\alpha \neq \beta$, we use \eqref{eq:cnt_Vi} to find out $(u,v), (v,u) \in E_G$ for any $u \in V_{i_\alpha}$ and $v\in V_{i_\beta}$. Thus, it follows that
\bq\label{eq:w_rel}
G[\mathcal{W}] = G[V_{i_1}] \iff \cdots \iff G[V_{i_\ell}] \in \mathcal{T}.
\eq
Now, we have three possibilities. If $\mathcal{W} =\emptyset$, then $G = G[\mathcal{U}]$ is semicomplete and hence $G \in \mathcal{T}$. If $\mathcal{U} = \emptyset$, then $G = G[\mathcal{W}] \in \mathcal{T}$ by \eqref{eq:w_rel}. Finally, if both $\mathcal{W}$ and $\mathcal{U}$ are nonempty, \eqref{eq:cnt_Vi} implies
\[
(u,w) \in E_G \quad \mbox{for all} \quad u \in \mathcal{U}, \ w \in \mathcal{W}.
\]
Therefore, $G \in \{ G[\mathcal{U}] \implies G[\mathcal{W}]\} \subseteq \mathcal{T}$ and this completes the induction, proving $\mathcal{S}\subseteq \mathcal{T}$.

\end{proof}

\bibliographystyle{abbrv}
\bibliography{sCS_network}

\end{document}